\documentclass[11pt]{article}
\usepackage{amsmath}
\usepackage{amssymb}

\usepackage{theorem}
\numberwithin{equation}{section}

\newtheorem{theorem}{Theorem}[section]

\newtheorem{corollary}[theorem]{Corollary}

\newtheorem{lemma}[theorem]{Lemma}
\newtheorem{definition}{Definition}[section]

\begin{document}
\title{Bilinear Strichartz estimates for the Schr{\"o}dinger map problem}
\date{\today}
\author{Benjamin Dodson}
\maketitle

\noindent \textbf{Abstract:} In this paper we prove bilinear Strichartz estimates for a solution to the Schr{\"o}dinger map problem whose size is small in the critical Strichartz space $\| |\nabla|^{\frac{d - 2}{2}} \psi_{x} \|_{L_{t,x}^{\frac{2(d + 2)}{d}}}$. These estimates will be useful in an upcoming paper in proving a local well - posedness result. Bilinear estimates make use of an argument similar to the argument found in \cite{PV}. We use the same gauges as in \cite{BIK}, \cite{BIKT}, and \cite{Smith}.

\section{Introduction}
The Schr{\"o}dinger map problem

\begin{equation}\label{1.1}
 \aligned
\partial_{t} \phi &= \phi \times \Delta_{x} \phi, \\
\phi(0) &= \phi_{0}, \\
\phi &: I \times \mathbf{R}^{d} \rightarrow S^{2} \hookrightarrow \mathbf{R}^{3}
\endaligned
\end{equation}

\noindent is a problem which has been a subject of a great deal of recent attention. This is a problem with a rich geometric structure that arises naturally in a number of different ways. See \cite{KLPST} or \cite{NSU} for more details.\vspace{5mm}

\noindent This system $(\ref{1.1})$ enjoys conservation of energy,

\begin{equation}\label{1.2}
 E(\phi(t)) = \frac{1}{2} \int_{\mathbf{R}^{d}} |\partial_{x} \phi(t,x)|^{2} dx
\end{equation}

\noindent and mass

\begin{equation}\label{1.3}
 M(\phi(t)) = \int_{\mathbf{R}^{d}} |\phi(t,x) - Q|^{2} dx,
\end{equation}

\noindent where $Q \in S^{2}$ is some fixed base point. When $d = 2$ both $(\ref{1.1})$ and $(\ref{1.2})$ are invariant with respect to the scaling

\begin{equation}\label{1.4}
 \phi(t,x) \mapsto \phi(\lambda^{2} t, \lambda x), \hspace{5mm} \lambda > 0.
\end{equation}

\noindent When $d = 2$ $(\ref{1.1})$ is called energy critical. \cite{BIK}, \cite{BIKT}, \cite{BIKT1}, \cite{Smith}, studied the partial differential equation satisfied by the derivatives of a solution to $(\ref{1.1})$. The derivatives of $\psi(t,x)$, $\psi_{l} = \partial_{x_{l}} \psi(t,x)$ satisfy an equation that is a perturbation of the free Schr{\"o}dinger equation

\begin{equation}\label{1.6}
 (i \partial_{t} + \Delta) \psi_{l} = -2i A_{m} \partial_{m} \psi_{l} -i (\partial_{m} A_{m}) \psi_{l} + (A_{t} + A_{m} A_{m}) \psi_{l} - i \psi_{m} Im(\bar{\psi}_{m} \psi_{l}).
\end{equation}

\noindent \textbf{Remark:} In this paper we adopt the usual convention that Latin letters $l, m = 1, ..., d$ and we sum over repeated indices. $A_{m}$, $A_{t}$ are the connection coefficients.\vspace{5mm}

\noindent Using the Coulomb gauge in dimensions $d \geq 4$ \cite{BIK} proved global well - posedness of $(\ref{1.1})$ for initial data sufficiently small in $\dot{H}^{d/2}$. \cite{BIKT} proved global well - posedness for small data in $d \geq 2$ using the caloric gauge. This result was subsequently extended by \cite{Smith} to data with energy below the energy of the ground state and $d = 2$, provided the data satisfies certain other smallness assumptions.\vspace{5mm}

\noindent The chief difficulty in the study of the derivative Schr{\"o}dinger maps equation arises from the magnetic term $A_{m} \partial_{m} \psi_{l}$ when $\psi_{l}$ is at a high frequency and $A_{m}$ is at a low frequency. This term cannot be treated perturbatively using only the Strichartz estimates. Instead \cite{BIK}, \cite{BIKT}, \cite{Smith} utilized bilinear Strichartz estimates to move half of the derivative from the high frequency term to the low frequency term. This combined with local smoothing results is enough to close the bootstrap under the smallness conditions of \cite{BIK}, \cite{BIKT}, and \cite{Smith}.\vspace{5mm}

\noindent In this paper we prove some bilinear Strichartz estimates for a solution to $(\ref{1.1})$. We start by recalling a bilinear Strichartz estimate for the linear Schr{\"o}dinger equation.

\begin{theorem}\label{t1.1}
 If $u$ solves the free Schr{\"o}dinger equation

\begin{equation}\label{1.7}
\aligned
i u_{t} + \Delta u = 0, \\
u(0) = u_{0},
\endaligned
\end{equation}

\noindent then for $M << N$, when $P_{N}$ is a Littlewood - Paley operator,

\begin{equation}\label{1.8}
 \| (P_{M} u)(P_{N} \bar{u}) \|_{L_{t,x}^{2}(\mathbf{R} \times \mathbf{R}^{d})} \lesssim \frac{M^{(d - 1)/2}}{N^{1/2}} \| P_{M} u_{0} \|_{L^{2}(\mathbf{R}^{d})} \| P_{N} u_{0} \|_{L^{2}(\mathbf{R}^{d})}.
\end{equation}

\end{theorem}

\noindent This can be proved using Fourier analytic techniques. \cite{B1} used the Fourier transform to prove this theorem when $d = 2$. The result was subsequently extended to all dimensions (see for example \cite{KilVis}). One can also prove a similar result on $I$ if $u$ solves

\begin{equation}\label{1.9}
i u_{t} + \Delta u = \pm |u|^{2} u,
\end{equation}

\begin{equation}\label{1.10}
 \| |\nabla|^{(d - 2)/2} u \|_{L_{t,x}^{\frac{2(d + 2)}{d}}(I \times \mathbf{R}^{d})} < \infty.
\end{equation}

\noindent \cite{PV} proved theorem $\ref{t1.1}$ via an interaction Morawetz estimate. This method is useful to this paper because it is very robust under perturbations of $(\ref{1.7})$. In particular, if $\psi_{l}$ solves $(\ref{1.6})$ then $(\ref{1.8})$ holds under a slight strengthening of $(\ref{1.10})$.\vspace{5mm}

\noindent First define a Sobolev space for $\phi : I \times \mathbf{R}^{d} \rightarrow S^{2}$.

\begin{definition}\label{d1.2}
Let $\mathcal F_{(d)}$ denote the Fourier transform on $L^{2}(\mathbf{R}^{d})$. For $\sigma \geq 0$ define the inhomogeneous Sobolev spaces on $\mathbf{R}^{d}$ for vector valued functions.

\begin{equation}\label{1.11}
H^{\sigma}(\mathbf{R}^{d}) = \{ f : \mathbf{R}^{d} \rightarrow \mathbf{C}^{n} : \| f \|_{H^{\sigma}(\mathbf{R}^{d})} = [\sum_{l = 1}^{n} \| \mathcal F(f_{l})(\xi) (1 + |\xi|^{2})^{\sigma/2} \|_{L^{2}(\mathbf{R}^{d})}^{2}]^{1/2} < \infty \},
\end{equation}

\noindent as well as the homogeneous Sobolev spaces

\begin{equation}\label{1.12}
\dot{H}^{\sigma}(\mathbf{R}^{d}) = \{ f : \mathbf{R}^{d} \rightarrow \mathbf{C}^{n} : \| f \|_{\dot{H}^{\sigma}(\mathbf{R}^{d})} = [\sum_{l = 1}^{n} \| \mathcal F(f_{l})(\xi) \cdot |\xi|^{\sigma} \|_{L^{2}(\mathbf{R}^{d})}^{2}]^{1/2} < \infty \}.
\end{equation}

\noindent For $\sigma \geq 0$, $Q = (Q_{1}, Q_{2}, Q_{3}) \in S^{2}$ define the complete metric space

\begin{equation}\label{1.13}
 H_{Q}^{\sigma}(\mathbf{R}^{d}; S^{2}) = \{ f : \mathbf{R}^{d} \rightarrow \mathbf{R}^{3} : |f(x)| \equiv 1, f - Q \in H^{\sigma} \}.
\end{equation}

\noindent This metric has the induced distance

\begin{equation}\label{1.14}
 d_{Q}^{\sigma}(f, g) = \| f - g \|_{H^{\sigma}(\mathbf{R}^{d})}.
\end{equation}

\noindent Let $\| f \|_{H_{Q}^{\sigma}} = d_{Q}^{\sigma}(f, Q)$ for $f \in H_{Q}^{\sigma}$. Define the complete metric spaces

\begin{equation}\label{1.15}
 H^{\infty} = H^{\infty}(\mathbf{R}^{d}; \mathbf{C}^{n}) = \cap_{\sigma \in \mathbf{Z}_{+}} H^{\sigma}(\mathbf{R}^{d}) \hspace{5mm} \text{and} \hspace{5mm} H_{Q}^{\infty}(\mathbf{R}^{d} ; S^{2}) = \cap_{\sigma \in \mathbf{Z}_{+}} H_{Q}^{\sigma}(\mathbf{R}^{d}; S^{2})
\end{equation}

\noindent with the induced distances.
\end{definition}

\noindent Choose a small constant $\delta > 0$, say $\delta = \frac{1}{40}$. Let $\psi_{x}$ be the vector $\psi_{x} = (\psi_{1}, ..., \psi_{d})$. Let $\beta(k)$ be a frequency envelope that majorizes $2^{k(d - 2)/2} \| P_{k} \psi_{x}(0) \|_{L_{x}^{2}(\mathbf{R}^{d})}$, satisfying

\begin{equation}\label{1.16}
2^{k(d - 2)/2} \| P_{k} \psi_{x}(0) \|_{L_{x}^{2}(\mathbf{R}^{d})} \leq \beta(k), \hspace{5mm} \beta(k) \leq 2^{\delta |k - l|} \beta(l), \hspace{5mm} \sum_{k} \beta(k)^{2} \lesssim \| \psi \|_{\dot{H}_{Q}^{d/2}}^{2}.
\end{equation}

\noindent For example one could choose

\begin{equation}\label{1.16.1}
 \beta(k) = \sum_{j} 2^{-\delta |j - k|} 2^{\frac{j(d - 2)}{2}} \| P_{j} \psi_{x}(0) \|_{L^{2}(\mathbf{R}^{d})}.
\end{equation}

\noindent Suppose also that $\alpha(k)$ is a frequency envelope that majorizes $2^{k(d - 2)/2} \| P_{k} \psi_{x} \|_{L_{t,x}^{\frac{2(d + 2)}{d}}(I \times \mathbf{R}^{d})}$, and

\begin{equation}\label{1.17}
2^{k(d - 2)/2} \| P_{k} \psi_{x} \|_{L_{t, x}^{\frac{2(d + 2)}{d}}(I \times \mathbf{R}^{d})} \leq \alpha(k), \hspace{5mm} \alpha(k) \leq 2^{\delta |k - l|} \alpha(l), \hspace{5mm} \sum_{k} \alpha(k)^{2} \lesssim \epsilon(\| \psi \|_{\dot{H}_{Q}^{d/2}})^{2}.
\end{equation}

\begin{theorem}\label{t1.3}
 Suppose $d \geq 4$ and $k - l \geq 10$, $\psi$ solves $(\ref{1.6})$, satisfies $(\ref{1.16})$ and $(\ref{1.17})$, and $A$ satisfies the Coulomb gauge. Then

\begin{equation}\label{1.18}
\| (P_{k} \bar{\psi}_{x})(P_{l} \psi_{x}) \|_{L_{t,x}^{2}(I \times \mathbf{R}^{d})} \lesssim 2^{-|l - k|/2} (\alpha(k) + \beta(k))(\alpha(l) + \beta(l)).
\end{equation}
\end{theorem}

\noindent \cite{BIK} was unable to use the Coulomb gauge in dimensions $d = 2$, $d = 3$. Instead, for dimensions $d \geq 2$ and small data \cite{BIKT} utilized the caloric gauge. \cite{Smith} also utilized the caloric gauge to study $d = 2$. The caloric gauge arises from computing the harmonic map heat flow with initial data $\psi_{x,t}(t,x)$ for any $(t,x) \in I \times \mathbf{R}^{d}$. The harmonic map heat flow is computed in for all $s > 0$. The gauge condition $A_{s} \equiv 0$ is imposed. Therefore it is necessary to study the bilinear estimates for $s, \tilde{s} \neq 0$ and $s \neq \tilde{s}$. \cite{Smith1} proved that the harmonic map heat flow is well - defined provided $\psi(t,x)$ has energy below the energy of the ground state. Since \cite{Smith1} only proved well - posedness of the harmonic map heat flow when $d = 2$,

\begin{theorem}\label{t1.4}
Suppose $d = 2$, $k - l \geq 10$, $\psi$ solves $(\ref{1.6})$, satisfies $(\ref{1.16})$ and $(\ref{1.17})$, and $A$ satisfies the caloric gauge. Moreover suppose $\psi(s,t,x)$ is the solution of the harmonic map heat flow with initial data $\psi(0,t,x)$. Then

\begin{equation}\label{1.19}
\| (P_{k} \bar{\psi}_{x}(s))(P_{l} \psi_{x}(\tilde{s})) \|_{L_{t,x}^{2}(I \times \mathbf{R}^{2})} \lesssim 2^{-|l - k|/2} (\alpha(k) + \beta(k))(\alpha(l) + \beta(l))(1 + s 2^{2k})^{-4} (1 + \tilde{s} 2^{2l})^{-4}.
\end{equation}
\end{theorem}

\noindent These results will be used in a subsequent paper to prove well - posedness of $(\ref{1.1})$.

\section{Gauge Field Equations}
\noindent Let $\phi$ be any function such that $\phi : \mathbf{R}^{2} \times (-T, T) \rightarrow S^{2}$. Denote space and time derivatives of $\phi$ as $\partial_{\alpha} \phi$, where $\alpha = 1, ..., d + 1$ and $\partial_{d + 1} \phi = \partial_{t} \phi$.\vspace{5mm}

\noindent \textbf{Remark:} The time variable is usually assigned to $\alpha = 0$. However this index will be reserved for time variable under the harmonic map heat flow in the caloric gauge.\vspace{5mm}

\noindent As in \cite{BIK}, \cite{BIKT}, and \cite{Smith} select an orthonormal frame $(v(t,x), w(t,x)) \in T_{\phi(t,x)} S^{2}$, i.e. smooth functions $v, w : \mathbf{R}^{2} \times (-T, T) \rightarrow S^{2}$ such that at each point $(x,t)$ the vectors $v(t,x)$, $w(t,x)$ form an orthonormal basis $T_{\phi(t,x)} S^{2}$. As a matter of convention assume $v$ and $w$ are chosen so that $v \times w = \phi$.\vspace{5mm}

\noindent Then introduce the derivative fields. Set

\begin{equation}\label{3.1}
 \psi_{\alpha} = v \cdot \partial_{\alpha} \phi + i w \cdot \partial_{\alpha} \phi.
\end{equation}

\noindent Then $\partial_{\alpha} \phi$ admits the representation

\begin{equation}\label{3.2}
 \partial_{\alpha} \phi = v Re(\psi_{\alpha}) + w Im(\psi_{\alpha}).
\end{equation}

\noindent Rewrite the vector $\partial_{\alpha} \phi$ with respect to the orthonormal basis $(v,w)$, then identify $\mathbf{R}^{2}$ with the complex numbers $\mathbf{C}$ according to $v \leftrightarrow 1$, $w \leftrightarrow i$. This identification respects the complex structure of the target manifold. The Riemannian connection on $S^{2}$ pulls back to a covariant derivative on $\mathbf{C}$, which we denote by

\begin{equation}\label{3.3}
 D_{\alpha} = \partial_{\alpha} + i A_{\alpha}.
\end{equation}

\noindent The connection coefficients $A_{\alpha}$ are defined via

\begin{equation}\label{3.4}
 A_{\alpha} = w \cdot \partial_{\alpha} v.
\end{equation}

\noindent Because the Riemannian connection on $S^{2}$ is torsion free the derivative fields satisfy the equations

\begin{equation}\label{3.5}
 D_{\beta} \psi_{\alpha} = D_{\alpha} \psi_{\beta}.
\end{equation}

\noindent Equivalently,

\begin{equation}\label{3.6}
 \partial_{\beta} A_{\alpha} - \partial_{\alpha} A_{\beta} = Im(\psi_{\beta} \bar{\psi}_{\alpha}) = q_{\beta \alpha}.
\end{equation}

\noindent If $\phi$ is a smooth solution to the Sch{\"o}dinger map problem $(\ref{1.1})$ then the derivatives satisfy the equation

\begin{equation}\label{3.7}
 \psi_{t} = i D_{l} \psi_{l}.
\end{equation}

\noindent This is because

\begin{equation}\label{3.8}
 \phi \times \Delta \phi = J(\phi) (\phi^{\ast} \nabla)_{j} \partial_{j} \phi,
\end{equation}

\noindent where $J(\phi)$ denotes the complex structure $\phi \times$ and $(\phi^{\ast} \nabla)_{j}$ the pullback of the Levi - Cevita connection $\nabla$ on the sphere. This implies

\begin{equation}\label{3.8.1}
\aligned
 (i \partial_{t} + \Delta) \psi_{l} &= -2i A_{m} \partial_{m} \psi_{l} -i (\partial_{m} A_{m}) \psi_{l} + (A_{t} + A_{m} A_{m}) \psi_{l} - i \psi_{m} Im(\bar{\psi}_{m} \psi_{l}), \\
D_{\alpha} \psi_{\beta} &= D_{\beta} \psi_{\alpha}, \\
Im(\psi_{\alpha} \bar{\psi}_{\beta}) &= \partial_{\alpha} A_{\beta} - \partial_{\beta} A_{\alpha}.
\endaligned
\end{equation}

\noindent A solution $\psi_{m}$ to $(\ref{3.7})$ cannot be determined uniquely without choosing an orthonormal frame $(v, w)$. Changing a given choice of orthonormal frame induces a gauge transformation and may be represented as

\begin{equation}\label{3.8.2}
 \psi_{m} \mapsto e^{i \theta} \psi_{m} \hspace{5mm} A_{m} \mapsto A_{m} + \partial_{m} \theta.
\end{equation}

\noindent The system $(\ref{3.7})$ is invariant with respect to such gauge transformations.\vspace{5mm}

\noindent In this paper we will discuss bilinear Strichartz estimates for two choices of gauge, the Coulomb gauge and the caloric gauge. The Coulomb gauge is a gauge which is quite useful in high dimensions (see \cite{BIK}) and in low dimensions when some additional symmetry is imposed on the problem(see \cite{GK} and \cite{BIKT1}). In this paper we will discuss the Coulomb gauge for dimensions $d \geq 4$.\vspace{5mm}

\noindent However, the Coulomb gauge becomes very difficult to use in low dimensions for a general Schr{\"o}dinger map problem. Therefore for dimension $d = 2$ we will consider the caloric gauge. This gauge was introduced in \cite{Taowm} to study wave maps in hyperbolic space. The series of papers \cite{Tao3}, \cite{Tao4}, \cite{Tao5}, \cite{Tao6}, \cite{Tao7} then used this gauge to establish global regularity of wave maps in hyperbolic space. \cite{Taosm} suggested that the caloric gauge would be a suitable gauge in which to study Schr{\"o}dinger maps. \cite{BIKT} utilized this gauge to establish global well - posedness in the setting of initial data with small critical norm. This result was further expanded by \cite{Smith}.\vspace{5mm}

\subsection{Coulomb Gauge:} Under the Coulomb gauge

\begin{equation}\label{3.9}
 \sum_{m = 1}^{d} \partial_{m} A_{m} = 0.
\end{equation}

\noindent In view of $(\ref{3.6})$ this leads to

\begin{equation}\label{3.10}
 A_{m} = \Delta^{-1} \sum_{l = 1}^{d} \partial_{l} Im(\bar{\psi}_{l} \psi_{m}).
\end{equation}

\noindent Also by $(\ref{3.6})$

\begin{equation}\label{3.11}
\Delta A_{0} = \sum_{l = 1}^{d} \partial_{l} (\partial_{0} A_{l} + Im(\psi_{l} \bar{\psi}_{d + 1})) = \sum_{l = 1}^{d} \partial_{l} Im(\psi_{l} \bar{\psi}_{d + 1}).
\end{equation}

\noindent Using $(\ref{3.5})$, $(\ref{3.7})$,

\begin{equation}\label{3.12}
 = -\sum_{m = 1}^{d} Re(\bar{\psi}_{l} D_{m} \psi_{m}) = -\sum_{m, l = 1}^{d} \partial_{l} \partial_{m} Re(\bar{\psi}_{l} \psi_{m}) + \frac{1}{2} \Delta (\sum_{m = 1}^{d} \psi_{m} \bar{\psi}_{m}).
\end{equation}

\noindent The caloric gauge will be discussed in an upcoming section.

\section{Proof of theorem $\ref{t1.1}$}
Everything in this section can be found in \cite{PV}. Theorem $\ref{t1.1}$ will be proved here for the reader's convenience, since the proof will be modified to deal with the case when $\psi$ solves $(\ref{1.6})$.\vspace{5mm}

\noindent Suppose $u$ solves

\begin{equation}\label{2.1}
(i \partial_{t} + \Delta) u = 0.
\end{equation}

\noindent The argument of \cite{PV} is more useful for this paper than the argument of \cite{B1} because it is very robust under perturbations of the Laplacian $\Delta$ or perturbations of $(\ref{2.1})$. Define the Morawetz potential

\begin{equation}\label{2.2}
\aligned
M(t) = \int |u_{M}(t,y)|^{2} \frac{(x - y)_{j}}{|x - y|} Im[\bar{u}_{N}(t,x) \partial_{j} u_{N}(t,x)] dx dy \\
+ \int |u_{N}(t,y)|^{2} \frac{(x - y)_{j}}{|x - y|} Im[\bar{u}_{M}(t,x) \partial_{j} u_{M}(t,x)] dx dy.
\endaligned
\end{equation}

\noindent Because $e^{it \Delta}$ is a Fourier multiplier and $|e^{-it |\xi|^{2}}| = 1$,

\begin{equation}\label{2.3}
\| u_{M}(t,x) \|_{L^{2}(\mathbf{R}^{d})} = \| u_{M}(0,x) \|_{L^{2}(\mathbf{R}^{d})},
\end{equation}

\noindent and therefore since $|\frac{(x - y)}{|x - y|}| \leq 1$,

\begin{equation}\label{2.4}
|M(t)| \lesssim (M + N) \| u_{M}(0,x) \|_{L^{2}(\mathbf{R}^{d})}^{2} \| u_{N}(0,x) \|_{L^{2}(\mathbf{R}^{d})}^{2}.
\end{equation}

\begin{lemma}\label{l2.1}
For $\omega \in S^{d - 1}$ let $x_{\omega} = x \cdot \omega$, $\partial_{\omega} = (\omega \cdot \nabla)$,

\begin{equation}\label{2.5}
\int_{S^{d - 1}} \frac{x_{\omega}}{|x_{\omega}|} f(x) \partial_{\omega} g(x) d\omega = \frac{1}{|x|} f(x) (x \cdot \nabla) g(x).
\end{equation}
\end{lemma}

\noindent \emph{Proof:} Without loss of generality suppose $x = (x_{1}, 0, ..., 0)$. 

\begin{equation}\label{2.6}
\frac{x_{\omega}}{|x_{\omega}|} = \frac{\omega_{1}}{|\omega_{1}|} \frac{x_{1}}{|x_{1}|}.
\end{equation}

\begin{equation}\label{2.7}
\int_{S^{d - 1}} \frac{x_{1}}{|x_{1}|} \frac{\omega_{1}}{|\omega_{1}|} f(x) \omega_{j} \partial_{j} f(x) = C(d) \frac{x_{1}}{|x_{1}|} f(x) \partial_{1} g(x) = C(d) \frac{x_{j}}{|x|} f(x) \partial_{j} g(x).
\end{equation}

\noindent Therefore, $M(t) = \frac{1}{C(d)} \int_{S^{d - 1}} M_{\omega}(t) d\omega$, where

\begin{equation}\label{2.8}
\aligned
M_{\omega}(t) = \int |u_{M}(t,y)|^{2} \frac{(x - y)_{\omega}}{|(x - y)_{\omega}|} Im[\bar{u}_{N}(t,x) \partial_{\omega} u_{N}(t,x)] dx dy \\
+ \int |u_{N}(t,y)|^{2}\frac{(x - y)_{\omega}}{|(x - y)_{\omega}|}  Im[\bar{u}_{M}(t,x) \partial_{\omega} u_{M}(t,x)] dx dy.
\endaligned
\end{equation}

\noindent By the fundamental theorem of calculus

\begin{equation}\label{2.9}
 M_{\omega}(T) - M_{\omega}(0) = \int_{0}^{T} \frac{d}{dt} M_{\omega}(t) dt.
\end{equation}

\noindent Without loss of generality take $\omega = (1, 0, ..., 0)$.

\begin{equation}\label{2.10}
\frac{d}{dt} M_{\omega}(t) = -2 \int \partial_{k} Im(\overline{u}_{M} \partial_{k} u_{M})(t,y) \frac{(x - y)_{1}}{|(x - y)_{1}|} Im[\bar{u}_{N} \partial_{1} u_{N}](t,x) dx dy
\end{equation}

\begin{equation}\label{2.11}
-2 \int \partial_{k} Im(\overline{u}_{N} \partial_{k} u_{N})(t,y) \frac{(x - y)_{1}}{|(x - y)_{1}|} Im[\bar{u}_{M} \partial_{1} u_{M}](t,x) dx dy
\end{equation}

\begin{equation}\label{2.12}
+ \frac{1}{2} \int |u_{M}(t,y)|^{2} \frac{(x - y)_{1}}{|(x - y)_{1}|} \partial_{1} \partial_{k}^{2} (|u_{N}(t,x)|^{2}) dx dy
\end{equation}

\begin{equation}\label{2.13}
+ \frac{1}{2} \int |u_{N}(t,y)|^{2} \frac{(x - y)_{1}}{|(x - y)_{1}|} \partial_{1} \partial_{k}^{2} (|u_{M}(t,x)|^{2}) dx dy
\end{equation}

\begin{equation}\label{2.14}
-2 \int |u_{M}(t,y)|^{2} \frac{(x - y)_{1}}{|(x - y)_{1}|} \partial_{k} Re(\partial_{1} \bar{u}_{N} \partial_{k} u_{N})(t,x) dx dy
\end{equation}

\begin{equation}\label{2.15}
-2 \int |u_{N}(t,y)|^{2} \frac{(x - y)_{1}}{|(x - y)_{1}|} \partial_{k} Re(\partial_{1} \bar{u}_{M} \partial_{k} u_{M})(t,x) dx dy.
\end{equation}

\noindent Integrating by parts

\begin{equation}\label{2.16}
= -2 \int  Im[\overline{u}_{M} \partial_{1} u_{M}](t,x_{1}, y_{2}, ..., y_{d})  Im[\bar{u}_{N} \partial_{1} u_{N}](t,x_{1}, x_{2}, ..., x_{d}) dx dy
\end{equation}

\begin{equation}\label{2.17}
-2 \int  Im[\overline{u}_{N} \partial_{1} u_{N}](t,x_{1}, y_{2}, ..., y_{d})  Im[\bar{u}_{M} \partial_{1} u_{M}](t,x_{1}, x_{2}, ..., x_{d}) dx dy
\end{equation}

\begin{equation}\label{2.18}
+ \frac{1}{2} \int \partial_{1}(|u_{M}(t,x_{1}, y_{2}, ..., y_{d})|^{2})  \partial_{1} (|u_{N}(t,x_{1}, x_{2}, ..., x_{d})|^{2}) dx dy
\end{equation}

\begin{equation}\label{2.19}
+ \frac{1}{2} \int \partial_{1}(|u_{N}(t,x_{1}, y_{2}, ..., y_{d})|^{2})  \partial_{1}  (|u_{M}(t,x_{1}, x_{2}, ..., x_{d})|^{2}) dx dy
\end{equation}

\begin{equation}\label{2.20}
+ 2 \int |u_{M}(t,x_{1}, y_{2}, ..., y_{d})|^{2}  |\partial_{1} u_{N}(t,x_{1}, x_{2}, ..., x_{d})|^{2} dx dy
\end{equation}

\begin{equation}\label{2.21}
+ 2 \int |u_{N}(t,x_{1}, y_{2}, ..., y_{d})|^{2}  |\partial_{1} u_{N}(t, x_{1}, x_{2}, ..., x_{d})|^{2} dx dy.
\end{equation}

\begin{equation}\label{2.22}
= \int \int |\partial_{1}( \bar{u}_{N}(t,x_{1}, x_{2}, ..., x_{d}) u_{M}(t, x_{1}, y_{2}, ..., y_{d})) |^{2} dx_{1} dx_{2} \cdots dx_{d} dy_{2} \cdots dy_{d}.
\end{equation}

\noindent In one dimension this implies

\begin{equation}\label{2.23}
\int \int |\partial_{x} (\bar{u}_{N} u_{M}) (t,x)|^{2} dx dt \lesssim (M + N) \| u_{M}(0) \|_{L^{2}(\mathbf{R})}^{2} \| u_{N}(0) \|_{L^{2}(\mathbf{R})}^{2}.
\end{equation}

\noindent Therefore Bernstein's inequality implies that when $M << N$,

\begin{equation}\label{2.24}
\| \bar{u}_{M} u_{N} \|_{L_{t,x}^{2}(\mathbf{R} \times \mathbf{R})} \lesssim \frac{1}{N^{1/2}} \| u_{M}(0) \|_{L^{2}(\mathbf{R})} \| u_{N}(0) \|_{L^{2}(\mathbf{R})},
\end{equation}

\noindent which concludes the proof of theorem $\ref{t1.1}$ when $d = 1$. In higher dimensions let $\tilde{P}_{M}$ be the Littlewood - Paley projection onto frequencies $|\xi_{2} + ... + \xi_{d}| \leq 100M$. This implies that for some $\phi(x)$, $|\phi(x)| \lesssim 1$, $|\phi(x)| \lesssim_{N} (1 + |x|)^{-N}$ for any $N$,

\begin{equation}\label{2.25}
 u_{M} = \tilde{P}_{M} u_{M} = \int_{\mathbf{R}^{d - 1}} u_{M}(x_{1}, x_{2} - \bar{y}) \phi(M \bar{y}) M^{d - 1} d\bar{y}.
\end{equation}

\begin{equation}\label{2.26}
 \partial_{1}(\bar{u}_{N}(t, x_{1}, x') u_{M}(t, x_{1}, x' + y_{0})) = \int_{\mathbf{R}^{d - 1}} \partial_{1}(\bar{u}_{N}(t, x_{1}, x') u_{M}(t, x_{1}, x' + y_{0} + \bar{y})) \phi(M \bar{y}) M^{d - 1} d\bar{y}.
\end{equation}

\noindent By Holder's inequality

\begin{equation}\label{2.27}
\aligned
 | \partial_{1}(\bar{u}_{N}(t, x_{1}, x') u_{M}(t, x_{1}, x' + y_{0})) | \\ \lesssim M^{\frac{d - 1}{2}} (\int_{\mathbf{R}^{d - 1}} |\partial_{1}(\bar{u}_{N}(t, x_{1}, x') u_{M}(t, x_{1}, x' + y_{0} + \bar{y}))|^{2} d\bar{y})^{1/2} (\int |\phi(M \bar{y})|^{2} M^{d - 1} d\bar{y})^{1/2}.
\endaligned
\end{equation}

\noindent Therefore by $(\ref{2.4})$, $(\ref{2.22})$,

\begin{equation}\label{2.28}
 \| \partial_{1}(\bar{u}_{N}(t, x_{1}, x') u_{M}(t, x_{1}, x' + y_{0})) \|_{L_{t,x}^{2}}^{2} \lesssim M^{d - 1} N \| u_{M}(0) \|_{L^{2}(\mathbf{R}^{d})}^{2} \| u_{N}(0) \|_{L^{2}(\mathbf{R}^{d})}^{2}.
\end{equation}

\noindent Integrating over $\omega \in S^{d - 1}$ implies

\begin{equation}\label{2.29}
 \| \nabla (\bar{u}_{N}(t, x_{1}, x') u_{M}(t, x_{1}, x' + y_{0})) \|_{L_{t,x}^{2}} \lesssim M^{\frac{d - 1}{2}} N^{1/2} \| u_{M}(0) \|_{L^{2}(\mathbf{R}^{d})} \| u_{N}(0) \|_{L^{2}(\mathbf{R}^{d})}.
\end{equation}

\noindent Applying Bernstein's inequality proves theorem $\ref{t1.1}$. $\Box$\vspace{5mm}

\noindent An identical computation would produce the same result with $u_{M}$ replaced by $u_{M}(x + x_{0})$ for some $x_{0} \in \mathbf{R}^{d}$. Therefore,

\begin{corollary}\label{c2.1}
 If $u$ solves the free Schr{\"o}dinger equation then for $M << N$,

\begin{equation}\label{2.30}
 \| (P_{M} u(t,x))(P_{N} \bar{u}(t,x + x_{0}) \|_{L_{t,x}^{2}(\mathbf{R} \times \mathbf{R}^{d})} \lesssim \frac{M^{(d - 1)/2}}{N^{1/2}} \| P_{M} u_{0} \|_{L^{2}(\mathbf{R}^{d})} \| P_{N} u_{0} \|_{L^{2}(\mathbf{R}^{d})}.
\end{equation}
\end{corollary}

\section{Almost Conserved Quantities}
Conservation of energy implies

\begin{equation}\label{4.14.1}
\| \psi_{x}(t) \|_{L_{x}^{2}(\mathbf{R}^{d})} = \| \psi_{x}(0) \|_{L_{x}^{2}(\mathbf{R}^{d})}.
\end{equation}

\noindent Therefore consider $d \geq 4$.

\begin{theorem}\label{t4.1}
For $d \geq 4$, $\epsilon(\| \psi \|_{\dot{H}_{Q}^{d/2}}) > 0$ sufficiently small,

\begin{equation}\label{4.14}
 \| |\nabla|^{\frac{d - 2}{2}} \psi_{x} \|_{L_{t}^{\infty} L_{x}^{2}(I \times \mathbf{R}^{d})} \lesssim \| |\nabla|^{\frac{d - 2}{2}} \psi_{x}(0) \|_{L^{2}(\mathbf{R}^{d})}.
\end{equation}
\end{theorem}

\noindent \emph{Proof:} Suppose $\psi_{x}$ solves $(\ref{1.6})$. Take the inner product

\begin{equation}\label{4.1}
 \langle u, v \rangle = Re \int u(x) \bar{v}(x) dx.
\end{equation}

\begin{equation}\label{4.2}
 \aligned
\frac{1}{2} \frac{d}{dt} \langle |\nabla|^{\frac{d - 2}{2}} \psi_{x}, |\nabla|^{\frac{d - 2}{2}} \psi_{x} \rangle = \langle i \Delta |\nabla|^{\frac{d - 2}{2}} \psi_{x}, |\nabla|^{\frac{d - 2}{2}} \psi_{x} \rangle \\
+ \langle |\nabla|^{\frac{d - 2}{2}} (-2 A_{m} \partial_{m} \psi_{l} - (\partial_{m} A_{m}) \psi_{l} - i (A_{t} + A_{m} A_{m}) \psi_{l} - \psi_{m} Im(\bar{\psi}_{m} \psi_{l})), |\nabla|^{\frac{d - 2}{2}} \psi_{l} \rangle.
\endaligned
\end{equation}

\noindent The first term on the right hand side of $(\ref{4.2})$ is $\equiv 0$.

\begin{equation}\label{4.3}
\| |\nabla|^{\frac{d - 2}{2}} ((-\partial_{m} A_{m}) \psi_{l} - i (A_{t} + A_{m} A_{m}) \psi_{l} - \psi_{m} Im(\bar{\psi}_{m} \psi_{l})) \|_{L_{t,x}^{\frac{2(d + 2)}{d + 4}}(I \times \mathbf{R}^{d})}
\end{equation}

\begin{equation}\label{4.4}
 \lesssim \| |\nabla \cdot A| + |A_{t}| + |A_{x} \cdot A_{x}| + |\psi_{x}|^{2} \|_{L_{t,x}^{\frac{d + 2}{2}}} \| |\nabla|^{\frac{d - 2}{2}} \psi_{x} \|_{L_{t,x}^{\frac{2(d + 2)}{d}}}
\end{equation}

\begin{equation}\label{4.5}
 + \| \psi_{x} \|_{L_{t,x}^{d + 2}} (\| |\nabla|^{\frac{d - 2}{2}} (\nabla \cdot A) \|_{L_{t,x}^{2}} + \| |\nabla|^{\frac{d - 2}{2}} A_{t} \|_{L_{t,x}^{2}} + \| ||\nabla|^{\frac{d - 2}{2}} (A_{x})||A_{x}| \|_{L_{t,x}^{2}})
\end{equation}

\begin{equation}\label{4.6}
 \lesssim \| |\nabla|^{\frac{d - 2}{2}} \psi_{x} \|_{L_{t,x}^{\frac{2(d + 2)}{d}}}^{\frac{d + 4}{d}} (\| |\nabla|^{\frac{d - 2}{2}} \psi_{x} \|_{L_{t}^{\infty} L_{x}^{2}}^{2 - 4/d} + \| |\nabla|^{\frac{d - 2}{2}} \psi_{x} \|_{L_{t}^{\infty} L_{x}^{2}}^{4 - 4/d}).
\end{equation}

\noindent Therefore by $(\ref{1.17})$,

\begin{equation}\label{4.7}
\int_{I} \langle |\nabla|^{\frac{d - 2}{2}} (- (\partial_{m} A_{m}) \psi_{l} - i (A_{t} + A_{m} A_{m}) \psi_{l} - \psi_{m} Im(\bar{\psi}_{m} \psi_{l})), |\nabla|^{\frac{d - 2}{2}} \psi_{l} \rangle dt
\end{equation}

\begin{equation}\label{4.8}
 \lesssim \epsilon^{\frac{2(d + 2)}{d}} (\| |\nabla|^{\frac{d - 2}{2}} \psi_{x} \|_{L_{t}^{\infty} L_{x}^{2}}^{2 - 4/d} + \| |\nabla|^{\frac{d - 2}{2}} \psi_{x} \|_{L_{t}^{\infty} L_{x}^{2}}^{4 - 4/d}).
\end{equation}

\noindent Finally evaluate

\begin{equation}\label{4.9}
 \int_{I} \langle 2 |\nabla|^{\frac{d - 2}{2}} A_{m} \partial_{m} \psi_{x}, |\nabla|^{\frac{d - 2}{2}} \psi_{x} \rangle dt = \int_{I} \langle 2 A_{m} \partial_{m} |\nabla|^{\frac{d - 2}{2}} \psi_{x}, |\nabla|^{\frac{d - 2}{2}} \psi_{x} \rangle dt
\end{equation}

\begin{equation}\label{4.10}
 + \int_{I} \langle 2 |\nabla|^{\frac{d - 2}{2}} A_{m} \partial_{m} \psi_{x}, |\nabla|^{\frac{d - 2}{2}} \psi_{x} \rangle - \langle 2 A_{m} \partial_{m} |\nabla|^{\frac{d - 2}{2}} \psi_{x}, |\nabla|^{\frac{d - 2}{2}} \psi_{x} \rangle dt.
\end{equation}

\noindent Integrate the right hand side of $(\ref{4.9})$ by parts.

\begin{equation}\label{4.11}
\aligned
 \int_{I} \int A_{m}(t,x) \partial_{m} ||\nabla|^{\frac{d - 2}{2}} \psi_{x}|^{2} dx dt = -\int_{I} \int (\nabla \cdot A) ||\nabla|^{\frac{d - 2}{2}} \psi_{x}|^{2} dx dt \\ \lesssim \| |\nabla|^{\frac{d - 2}{2}} \psi_{x} \|_{L_{t,x}^{\frac{2(d + 2)}{d}}}^{\frac{2(d + 2)}{d}} \| |\nabla|^{\frac{d - 2}{2}} \|_{L_{t}^{\infty} L_{x}^{2}}^{2 - 4/d}.
\endaligned
\end{equation}

\noindent Therefore,

\begin{equation}\label{4.12}
 (\ref{4.10}) \lesssim \| (|\nabla|^{\frac{d - 2}{2}} A_{m}) (\partial_{m} \psi_{x}) \|_{L_{t,x}^{\frac{2(d + 2)}{d + 4}}} \| |\nabla|^{\frac{d - 2}{2}} \psi_{x} \|_{L_{t,x}^{\frac{2(d + 2)}{d}}} + \| \nabla A_{m} \|_{L_{t,x}^{\frac{d + 2}{2}}} \| |\nabla|^{\frac{d - 2}{2}} \psi_{x} \|_{L_{t,x}^{\frac{2(d + 2)}{d}}}^{2}
\end{equation}

\begin{equation}\label{4.13}
 \lesssim \| |\nabla|^{\frac{d - 2}{2}} \psi_{x} \|_{L_{t,x}^{\frac{2(d + 2)}{d}}}^{\frac{2(d + 2)}{d}} \| |\nabla|^{\frac{2(d + 2)}{d}} \psi_{x} \|_{L_{t}^{\infty} L_{x}^{2}}^{2 - 4/d}.
\end{equation}

\noindent Putting together $(\ref{4.8})$, $(\ref{4.11})$, and $(\ref{4.13})$, for $\epsilon(\| |\nabla|^{\frac{d - 2}{2}} \psi_{x}(0) \|_{L^{2}}) > 0$ sufficiently small, say

\begin{equation}\label{4.15}
 \epsilon^{4/d} \| |\nabla|^{\frac{d - 2}{2}} \psi_{x}(0) \|_{L^{2}}^{4 - 4/d} << 1,
\end{equation}

\noindent the theorem is proved. $\Box$\vspace{5mm}

\noindent In order to make use of the interaction Morawetz estimates of the previous section we need to estimate $\| P_{N} \psi_{x}(t) \|_{L_{x}^{2}(\mathbf{R}^{d})}$ when $t \in I$.

\begin{equation}\label{4.16}
 \frac{1}{2} \frac{d}{dt} \langle P_{M} \psi_{x}, P_{M} \psi_{x} \rangle = -\langle 2 P_{M}(A_{m} \partial_{m} \psi_{x}), P_{M} \psi_{x} \rangle
\end{equation}

\begin{equation}\label{4.17}
 - \langle P_{M}((\partial_{m} A_{m}) \psi_{x}), P_{M} \psi_{x} \rangle - \langle i P_{M}((A_{t} + A_{x} \cdot A_{x}) \psi_{x}), P_{M} \psi_{x} \rangle - \langle P_{M} (\psi_{m} Im(\bar{\psi}_{m} \psi_{l})), P_{M} \psi_{l} \rangle.
\end{equation}

\begin{lemma}\label{l4.2}
 When $d \geq 4$, $\psi$ solves $(\ref{1.6})$ and satisfies $(\ref{1.16})$ and $(\ref{1.17})$, $M = 2^{k}$ for some $k \in \mathbf{Z}$,

\begin{equation}\label{4.18}
 \| P_{M}((P_{\geq \frac{M}{100}} A_{m}) \partial_{m} \psi_{x} + (\partial_{m} A_{m}) \psi_{x} + i(A_{t} + A_{x} \cdot A_{x}) \psi_{x} + \psi_{m} Im(\bar{\psi}_{m} \psi_{x})) \|_{L_{t,x}^{\frac{2(d + 2)}{d + 4}}} \lesssim 2^{-k(d - 2)/2} \alpha(k).
\end{equation}

\end{lemma}

\noindent \textbf{Remark:} If $M = 2^{k}$, $N = 2^{j}$, let $\alpha(M) = \alpha(k)$ and $\alpha(N) = \alpha(j)$.\vspace{5mm}

\noindent \emph{Proof:} Begin with the easiest term and move to the most difficult. By $(\ref{1.17})$

\begin{equation}\label{4.19}
 \| P_{M}(\psi_{m} Im(\bar{\psi}_{m} \psi_{x})) \|_{L_{t,x}^{\frac{2(d + 2)}{d + 4}}} \lesssim \| \psi_{x} \|_{L_{t,x}^{\frac{d + 2}{2}}} \| P_{\geq \frac{M}{100}} \psi_{x} \|_{L_{t,x}^{\frac{2(d + 2)}{d}}} \lesssim \epsilon^{4/d} \| |\nabla|^{\frac{d - 2}{2}} \psi_{x} \|_{L_{t}^{\infty} L_{x}^{2}}^{2 - 4/d} \| P_{\geq \frac{M}{100}} \psi_{x} \|_{L_{t,x}^{\frac{2(d + 2)}{d}}}.
\end{equation}

\noindent By $(\ref{4.15})$

\begin{equation}\label{4.20}
 \lesssim \sum_{l \geq k - 10}^{\infty} 2^{-l(d - 2)/2} \alpha(l) \lesssim 2^{-k(d - 2)/2} \alpha(k).
\end{equation}

\noindent Likewise,

\begin{equation}\label{4.21}
\| P_{M}((\partial_{m} A_{m}) \psi_{x} + i(A_{t} + A_{x} \cdot A_{x}) \psi_{x}) \|_{L_{t,x}^{\frac{2(d + 2)}{d + 4}}} \lesssim 2^{-k(d - 2)/2} \alpha(k)
\end{equation}

\begin{equation}\label{4.22}
 + \| P_{\geq \frac{M}{4}} ((\partial_{m} A_{m}) + i(A_{t} + A_{x} \cdot A_{x})) \|_{L_{t,x}^{2}} \| P_{\leq \frac{M}{100}} \psi_{x} \|_{L_{t,x}^{d + 2}}.
\end{equation}

\noindent By $(\ref{3.10})$, $(\ref{3.12})$, and $(\ref{4.20})$,

\begin{equation}\label{4.23}
\aligned
 \| P_{\geq \frac{M}{4}} ((\partial_{m} A_{m}) + i(A_{t} + A_{x} \cdot A_{x})) \|_{L_{t,x}^{2}} \\ \lesssim \| P_{\geq \frac{M}{100}} \psi_{x} \|_{L_{t,x}^{\frac{2(d + 2)}{d}}} \| \psi_{x} \|_{L_{t,x}^{\frac{2(d + 2)}{d}}}^{4/d} (\| |\nabla|^{\frac{d - 2}{2}} \psi_{x} \|_{L_{t}^{\infty} L_{x}^{2}}^{2 - 4/d} + \| |\nabla|^{\frac{d - 2}{2}} \psi_{x} \|_{L_{t}^{\infty} L_{x}^{2}}^{4 - 4/d}) \lesssim 2^{-k(d - 2)/2} \alpha(k).
\endaligned
\end{equation}

\noindent $\Box$\vspace{5mm}

\noindent Now consider

\begin{equation}\label{4.24}
 \int_{I} \langle 2 P_{M} ((P_{\leq \frac{M}{100}} A_{m}) \partial_{m} \psi_{x}), P_{M} \psi_{x} \rangle dt = 2 \int_{I} \langle (P_{\leq \frac{M}{100}} A_{m}) \partial_{m} (P_{M} \psi_{x}), P_{M} \psi_{x} \rangle dt
\end{equation}

\begin{equation}\label{4.25}
 + 2 \int_{I} \langle P_{M} ((P_{\leq \frac{M}{100}} A_{m}) \partial_{m} \psi_{x}), P_{M} \psi_{x} \rangle dt - 2 \int_{I} \langle (P_{\leq \frac{M}{100}} A_{m}) \partial_{m} (P_{M} \psi_{x}), P_{M} \psi_{x} \rangle dt.
\end{equation}

\noindent Integrating the right hand side of $(\ref{4.24})$ by parts

\begin{equation}\label{4.26}
\aligned
 \int_{I} \int (P_{\leq \frac{M}{100}} A_{m}) \partial_{m} |P_{M} \psi_{x}|^{2} dx dt \lesssim \| \nabla \cdot A \|_{L_{t,x}^{\frac{d + 2}{2}}} \| P_{M} \psi_{x} \|_{L_{t,x}^{\frac{2(d + 2)}{d}}}^{2} \\ \lesssim \epsilon^{4/d} \| |\nabla|^{\frac{d - 2}{2}} \psi_{x} \|_{L_{t}^{\infty} L_{x}^{2}}^{2 - 4/d} \alpha(k)^{2} 2^{-k(d - 2)} \lesssim \alpha(k)^{2} 2^{-k(d - 2)}.
\endaligned
\end{equation}

\noindent By the fundamental theorem of calculus we have the estimate on the Littlewood - Paley multiplier for $|\eta| << M$, $|\xi| \sim M$.

\begin{equation}\label{4.27}
 |\phi(\frac{\xi}{M}) - \phi(\frac{\xi + \eta}{M})| \lesssim \frac{1}{M} |\eta|.
\end{equation}

\noindent Therefore,

\begin{equation}\label{4.28}
 \| P_{M}((P_{\leq \frac{M}{100}} A_{m}) \partial_{m} \psi_{x}) - (P_{\leq \frac{M}{100}} A_{m}) \partial_{m} P_{M} \psi_{x} \|_{L_{t,x}^{\frac{2(d + 2)}{d + 4}}} \lesssim \| P_{M} \psi_{x} \|_{L_{t,x}^{\frac{2(d + 2)}{d}}} \| \partial_{x} A \|_{L_{t,x}^{\frac{d + 2}{d}}} \lesssim \alpha(k) 2^{-k(d - 2)/2}.
\end{equation}

\noindent Combining lemma $\ref{l4.2}$, $(\ref{4.26})$, and $(\ref{4.28})$

\begin{theorem}\label{t4.3}
 If $\psi_{x}$ satisfies $(\ref{1.16})$, $(\ref{1.17})$, and $(\ref{4.15})$ also holds, for $t \in I$,

\begin{equation}\label{4.29}
 \| P_{M} \psi_{x}(t) \|_{L_{x}^{2}} \lesssim M^{-\frac{d - 2}{2}} (\alpha(M) + \beta(M)).
\end{equation}

\end{theorem}

\section{Proof of theorem $\ref{t1.3}$:}
Now suppose $\psi_{x}$ solves $(\ref{1.6})$ and $d \geq 4$. Theorem $\ref{t4.3}$, $(\ref{1.16})$, and $(\ref{1.17})$ imply

\begin{equation}\label{5.1}
 \sup_{t \in I} |M(t)| \lesssim N(\alpha(M) + \beta(M))^{2} (\alpha(N) + \beta(N))^{2}.
\end{equation}

\noindent By corollary $\ref{c2.1}$ this would automatically imply theorem $\ref{t1.3}$ if

\begin{equation}\label{5.2}
 (\partial_{t} - i \Delta) \psi_{x} = 0.
\end{equation}

\noindent Therefore it suffices to bound the errors arising from the right hand side of $(\ref{1.6})$. These errors are quite similar to the errors in the proof of theorem $\ref{t4.3}$. Without loss of generality it suffices to consider the error terms in

\begin{equation}\label{5.3}
 \int \int |P_{M} \psi_{x}(t,y)|^{2} \frac{(x - y)_{j}}{|x - y|} Im[P_{N} \bar{\psi}_{x}(t,x) \partial_{j} \psi_{x}(t,x)] dx dy.
\end{equation}

\noindent The error is given by

\begin{equation}\label{5.4}
\mathcal E = 2 \int \int \int Re[P_{M}(\overline{\psi_{x}}(t,y)) (\partial_{t} - i \Delta) P_{M}(\psi_{x}(t,y))] \frac{(x - y)_{j}}{|x - y|} Im[(P_{N} \bar{\psi}_{x}(t,x) \partial_{j} P_{N} \psi_{x}(t,x)] dx dy dt
\end{equation}

\begin{equation}\label{5.5}
 + \int \int \int |P_{M} \psi_{x}(t,y)|^{2} \frac{(x - y)_{j}}{|x - y|} Im[(P_{N} \bar{\psi}_{x}(t,x) \partial_{j} P_{N}(\partial_{t} - i \Delta) \psi_{x}(t,x)] dx dy dt
\end{equation}

\begin{equation}\label{5.6}
 + \int \int \int |P_{M} \psi_{x}(t,y)|^{2} \frac{(x - y)_{j}}{|x - y|} Im[(P_{N} \overline{(\partial_{t} - i\Delta) \psi_{x}(t,x)} \partial_{j} P_{N} \psi_{x}(t,x)] dx dy dt.
\end{equation}

\noindent By lemma $\ref{l4.2}$, since

\begin{equation}\label{5.7}
\| P_{N} (-2 (P_{\geq \frac{N}{100}} A_{m}) \partial_{m} \psi_{x} - (\partial_{m} A_{m}) \psi_{x} - i(A_{t} + A_{x} \cdot A_{x}) \psi_{x}) \|_{L_{t,x}^{\frac{2(d + 2)}{d + 4}}} \lesssim N^{-\frac{d - 2}{2}} \alpha(N),
\end{equation}

\begin{equation}\label{5.8}
\aligned
\mathcal E = -4 \int \int \int Re[P_{M}(\overline{\psi_{x}}(t,y))  P_{M}((P_{\leq \frac{M}{100}} A_{m}) \partial_{m} \psi_{x}(t,y))] \\ \frac{(x - y)_{j}}{|x - y|} Im[(P_{N} \bar{\psi}_{x}(t,x) \partial_{j} P_{N} \psi_{x}(t,x)] dx dy dt
\endaligned
\end{equation}

\begin{equation}\label{5.9}
 -2 \int \int \int |P_{M} \psi_{x}(t,y)|^{2} \frac{(x - y)_{j}}{|x - y|} Im[(P_{N} \bar{\psi}_{x}(t,x) \partial_{j} P_{N} (P_{\leq \frac{N}{100}} A_{m}) \partial_{m} \psi_{x}(t,x))] dx dy dt
\end{equation}

\begin{equation}\label{5.10}
 -2 \int \int \int |P_{M} \psi_{x}(t,y)|^{2} \frac{(x - y)_{j}}{|x - y|} Im[(P_{N} \overline{(P_{\leq \frac{N}{100} \partial_{m} A_{m})} \psi_{x}(t,x)} \partial_{j} P_{N} \psi_{x}(t,x)] dx dy dt
\end{equation}

\begin{equation}\label{5.11}
 + O(N N^{-(d - 2)} M^{-(d - 2)}) (\alpha(M) + \beta(M))^{2} (\alpha(N) + \beta(N))^{2}.
\end{equation}

\noindent Integrate $(\ref{5.8})$ and $(\ref{5.9}) + (\ref{5.10})$ by parts.

\begin{equation}\label{5.12}
 \int_{I} \int \int (P_{\leq \frac{M}{100}} A_{m}) \partial_{m} |P_{M} \psi_{x}|^{2} \frac{(x - y)_{j}}{|x - y|} Im[P_{N} \bar{\psi}_{x} \partial_{j} P_{N} \psi_{x}] dx dy dt
\end{equation}

\begin{equation}\label{5.13}
 = -\int_{I} \int \int \partial_{m} (P_{\leq \frac{M}{100}} A_{m}) |P_{M} \psi_{x}|^{2} \frac{(x - y)_{j}}{|x - y|} Im[P_{N} \bar{\psi}_{x} \partial_{j} P_{N} \psi_{x}] dx dy dt
\end{equation}

\begin{equation}\label{5.14}
 -\int_{I} \int \int  (P_{\leq \frac{M}{100}} A_{m}) |P_{M} \psi_{x}|^{2} \partial_{m}(\frac{(x - y)_{j}}{|x - y|}) Im[P_{N} \bar{\psi}_{x} \partial_{j} P_{N} \psi_{x}] dx dy dt.
\end{equation}

\begin{equation}\label{5.15}
\aligned
 (\ref{5.13}) \lesssim N \| P_{N} \psi_{x} \|_{L_{t}^{\infty} L_{x}^{2}}^{2} \| \partial_{x} A \|_{L_{t,x}^{\frac{d + 2}{2}}} \| P_{M} \psi_{x} \|_{L_{t,x}^{\frac{2(d + 2)}{d}}} \\ \lesssim N N^{-(d - 2)} M^{-(d - 2)} (\alpha(M) + \beta(M))^{2} (\alpha(N) + \beta(N))^{2}.
\endaligned
\end{equation}

\noindent The Hardy - Littlewood - Sobolev inequality implies

\begin{equation}\label{5.16}
 (\ref{5.14}) \lesssim \| P_{M} \psi_{x} \|_{L_{t,x}^{\frac{2(d + 2)}{d}}}^{2} \| P_{\leq \frac{M}{100}} A_{m} \|_{L_{t}^{\infty} L_{x}^{d}} \| P_{N} \psi_{x} \|_{L_{t,x}^{\frac{2(d + 2)}{d}}}^{4/d} \| P_{N} \psi_{x} \|_{L_{t}^{\infty} L_{x}^{2}}^{2 - 4/d}.
\end{equation}

\noindent Likewise

\begin{equation}\label{5.17}
 \int_{I} \int \int |P_{M} \psi_{x}(t,y)|^{2} \frac{(x - y)_{j}}{|x - y|} Im[(P_{N} \bar{\psi}_{x}) (P_{\leq \frac{N}{100}} A_{m}) \partial_{j} \partial_{m} P_{N} \psi_{x}](t,x) dx dy
\end{equation}

\begin{equation}\label{5.18}
 + \int_{I} \int \int |P_{M} \psi_{x}(t,y)|^{2} \frac{(x - y)_{j}}{|x - y|} Im[(P_{\leq \frac{N}{100}} A_{m}) \partial_{m} (P_{N} \bar{\psi}_{x})  \partial_{j}  P_{N} \psi_{x}](t,x) dx dy
\end{equation}

\begin{equation}\label{5.19}
 = -\int_{I} \int \int |P_{M} \psi_{x}(t,y)|^{2} \frac{(x - y)_{j}}{|x - y|} \partial_{m} (P_{\leq \frac{N}{100}} A_{m}) Im[(P_{N} \bar{\psi}_{x}) \partial_{j} P_{N} \psi_{x}](t,x) dx dy dt
\end{equation}

\begin{equation}\label{5.20}
 = -\int_{I} \int \int |P_{M} \psi_{x}(t,y)|^{2} \partial_{m}(\frac{(x - y)_{j}}{|x - y|})  (P_{\leq \frac{N}{100}} A_{m}) Im[(P_{N} \bar{\psi}_{x}) \partial_{j} P_{N} \psi_{x}](t,x) dx dy dt.
\end{equation}

\noindent Once again use the Hardy - Littlewood - Sobolev theorem for $(\ref{5.20})$. As in the proof of theorem $\ref{t4.3}$

\begin{equation}\label{5.21}
 \| P_{M}((P_{\leq \frac{M}{100}} A_{m}) \partial_{m} \psi_{x}) (P_{M} \bar{\psi}_{x}) - (P_{\leq \frac{M}{100}} A_{m})(P_{M} \partial_{m} \psi_{x})(P_{M} \bar{\psi}_{x}) \|_{L_{t,x}^{1}} \lesssim M^{-(d - 2)} \alpha(M)^{2}.
\end{equation}

\begin{equation}\label{5.22}
 \| (P_{N} \bar{\psi}_{x}) \partial_{j}((P_{\leq \frac{N}{100}} A_{m}) \partial_{m} (P_{N} \psi_{x})) - (P_{N} \bar{\psi}_{x})(P_{\leq \frac{N}{100}} A_{m}) \partial_{j} \partial_{m} (P_{N} \psi_{x}) \|_{L_{t,x}^{1}} \lesssim N \cdot N^{-(d - 2)} \alpha(N)^{2}.
\end{equation}

\begin{equation}\label{5.23}
 \| (P_{N} \bar{\psi}_{x}) \partial_{j}((P_{\leq \frac{N}{100}} A_{m}) \partial_{m} (P_{N} \psi_{x}) - P_{N}((P_{\leq \frac{N}{100}} A_{m}) \partial_{m}  \psi_{x})) \|_{L_{t,x}^{1}} \lesssim N \cdot N^{-(d - 2)} \alpha(N)^{2}.
\end{equation}

\begin{equation}\label{5.24}
 \|  (((P_{\leq \frac{N}{100}} A_{m}) \partial_{m} (P_{N} \psi_{x}) - P_{N}( (P_{\leq \frac{N}{100}} A_{m}) \partial_{m}  \psi_{x})) \partial_{j}(P_{N} \psi_{x})\|_{L_{t,x}^{1}} \lesssim N \cdot N^{-(d - 2)} \alpha(N)^{2}.
\end{equation}

\noindent This proves theorem $\ref{t1.3}$. $\Box$

\section{Caloric Gauge}
The caloric gauge was proposed in \cite{Taowm} in the context of wave maps and then in \cite{Taosm} in the context of Schr{\"o}dinger maps. Precisely, at each time $t$ we solve the covariant heat equation with $\phi(t)$ as the initial data on $[0, \infty) \times \mathbf{R}^{d}$,

\begin{equation}\label{6.0}
 \aligned
\partial_{s} \tilde{\phi} &= \Delta_{x} \tilde{\phi} + \tilde{\phi} \cdot \sum_{m = 1}^{d} |\partial_{m} \tilde{\phi}|^{2} \\
\tilde{\phi}(0, t,x) &= \phi(t,x).
\endaligned
\end{equation}

\noindent \cite{Smith1} proved that $(\ref{6.0})$ is well - posed on $\mathbf{R}^{2}$ for $s > 0$ when the energy of $\phi$ is less than the energy of the ground state. Moreover, $\phi(s)$ approaches the equilibrium state $Q$ as $s \rightarrow \infty$. Therefore we can choose $(v_{\infty}, w_{\infty})$ at $s = \infty$ as an arbitrary orthonormal base in $T_{Q} S^{2}$. Pulling back $(v_{\infty}, w_{\infty})$ along the backward heat flow by parallel transport gives an orthonormal frame $(v, w)$ for all $s \geq 0$. Moreover,

\begin{equation}\label{6.0.1}
w \cdot \partial_{s} v = A_{s} = 0.
\end{equation}

\noindent In the gauge the harmonic map heat flow is given by

\begin{equation}\label{6.1}
 (\partial_{s} - \Delta_{x}) \psi_{m} = 2i A_{l} \partial_{l} \psi_{m} - (A_{l} \cdot A_{l} - i \partial_{l} A_{l}) \psi_{m} + i Im(\psi_{m} \bar{\psi}_{l}) \psi_{l}.
\end{equation}

\begin{equation}\label{6.2}
 (\partial_{s} - \Delta_{x}) \psi_{t} = 2i A_{l} \partial_{l} \psi_{t} - (A_{l} \cdot A_{l} - i \partial_{l} A_{l}) \psi_{t} + i Im(\psi_{t} \bar{\psi}_{l}) \psi_{l}.
\end{equation}

\noindent $(\ref{3.6})$ implies

\begin{equation}\label{6.3}
 \partial_{0} A_{s} = Im(\psi_{0} \bar{\psi}_{m}).
\end{equation}

\noindent Integrating backward from $s = \infty$, for any $m = 1, ..., d + 1$,

\begin{equation}\label{6.4}
 A_{m}(s) = -\int_{s}^{\infty} Im(\bar{\psi}_{m} (\partial_{l} \psi_{l} + i A_{l} \psi_{l}))(r) dr.
\end{equation}

\begin{theorem}\label{t6.1}
Let $\phi$ be a heat flow with classical initial data whose energy $E_{0}$ is less than $E_{crit}$. Let $e$ be a caloric gauge for $\phi$, and let $A_{x}$ denote the connection fields. Then we have the pointwise bounds

\begin{equation}\label{6.5}
\sup_{s > 0} s^{(k + 1)/2} \| \partial_{x}^{k} A_{x}(s) \|_{L_{x}^{\infty}(\mathbf{R}^{2})} \lesssim_{E_{0}, k} 1,
\end{equation}

\begin{equation}\label{6.6}
 \sup_{s > 0} s^{k/2} \| \partial_{x}^{k} A_{x}(s) \|_{L_{x}^{2}(\mathbf{R}^{2})} \lesssim_{E_{0},k} 1,
\end{equation}

\noindent for all $k \geq 0$ and $s > 0$, as well as

\begin{equation}\label{6.7}
 \int_{0}^{\infty} s^{(k - 1)/2} \| \partial_{x}^{k} A_{x}(s) \|_{L_{x}^{\infty}(\mathbf{R}^{2})} ds \lesssim_{E_{0}, k} 1,
\end{equation}

\begin{equation}\label{6.8}
 \int_{0}^{\infty} s^{(k - 1)/2} \| \partial_{x}^{k + 1} A_{x}(s) \|_{L_{x}^{2}(\mathbf{R}^{2})} ds \lesssim_{E_{0},k} 1.
\end{equation}

\noindent For all $k \geq 0$.

\end{theorem}

\noindent \emph{Proof:} This was proved in theorem $7.4$ of \cite{Smith1}. $\Box$

\begin{corollary}\label{c6.2}
Let $\phi$ be a heat flow with classical initial data with energy $E_{0}$ less than $E_{crit}$. Let $e$ be a caloric gauge for $\phi$. Then for all $k \geq 1$,

\begin{equation}\label{6.9}
\int_{0}^{\infty} s^{k - 1} \| \partial_{x}^{k} \psi_{x} \|_{L_{x}^{2}(\mathbf{R}^{2})}^{2} \lesssim_{E_{0}, k} 1,
\end{equation}

\begin{equation}\label{6.10}
 \sup_{s > 0} s^{(k - 1)/2} \| \partial_{x}^{k - 1} \psi_{x} \|_{L_{x}^{2}(\mathbf{R}^{2})} \lesssim_{E_{0}, k} 1,
\end{equation}

\begin{equation}\label{6.11}
 \int_{0}^{\infty} s^{k - 1} \| \partial_{x}^{k - 1} \psi_{x} \|_{L_{x}^{\infty}(\mathbf{R}^{2})}^{2} ds \lesssim_{E_{0}, k} 1,
\end{equation}

\begin{equation}\label{6.12}
 \sup_{s > 0} s^{k/2} \| \partial_{x}^{k - 1} \psi_{x} \|_{L_{x}^{\infty}(\mathbf{R}^{2})} \lesssim_{E_{0}, k} 1.
\end{equation}

\noindent Analogous estimates hold if one replaces $\partial_{x} \psi_{x}$ with $\psi_{s}$, $\partial_{x}^{2}$ with $\partial_{s}$, and/or $\partial_{x}$ with $D_{x}$.

\end{corollary}

\noindent \emph{Proof:} This is corollary $7.5$ in \cite{Smith1}. $\Box$

\begin{theorem}\label{t6.3}
 For $t \in I$, $\alpha$ and $\beta$ satisfy $(\ref{1.16})$ and $(\ref{1.17})$, $d = 2$,

\begin{equation}\label{6.13}
 \| P_{M} \psi_{x}(t) \|_{L^{2}(\mathbf{R}^{2})}^{2} \leq \alpha(M)^{2} + C(E_{0})(\beta(M)^{2} \epsilon^{2} + \| (P_{M} \bar{\psi}_{x}) \psi_{x} \|_{L_{t,x}^{2}} \epsilon^{2}).
\end{equation}

\end{theorem}

\noindent \emph{Proof:}

\begin{equation}\label{6.14}
 \frac{d}{dt} \| P_{M} \psi_{x} \|_{L^{2}}^{2} = -\langle 2 P_{M}(A_{m} \partial_{m} \psi_{x}), P_{M} \psi_{x} \rangle
\end{equation}

\begin{equation}\label{6.15}
- \langle P_{M}((\partial_{m} A_{m}) \psi_{x}), P_{M} \psi_{x} \rangle - \langle i P_{M}((A_{t} + A_{x} \cdot A_{x}) \psi_{x}), P_{M} \psi_{x} \rangle - \langle P_{M}(\psi_{m} Im(\bar{\psi}_{m} \psi_{x})), P_{M} \psi_{x} \rangle.
\end{equation}

\noindent This implies that for $t \in I$,

\begin{equation}\label{6.16}
 \| P_{M} \psi_{x}(t) \|_{L^{2}(\mathbf{R}^{2})}^{2} \leq \beta(M)^{2} - \int_{I} \langle 2 P_{M} (A_{m} \partial_{m} \psi_{x}), P_{M} \psi_{x} \rangle dt
\end{equation}

\begin{equation}\label{6.17}
 + \| (P_{M} \bar{\psi}_{x}) \psi_{x} \|_{L_{t,x}^{2}} (\| \partial_{x} A \|_{L_{t,x}^{2}} + \| A_{t} \|_{L_{t,x}^{2}} + \| A_{x} \|_{L_{t,x}^{4}}^{2} + \| \psi_{x} \|_{L_{t,x}^{4}}^{2}).
\end{equation}

\noindent As in the Coulomb gauge

\begin{equation}\label{6.20}
 -2 \int_{I} \langle P_{M}((P_{\leq \frac{M}{100}} A_{m}) \partial_{m} \psi_{x}), P_{M} \psi_{x} \rangle dt + 2 \int_{I} \langle (P_{\leq \frac{M}{100}} A_{m}) \partial_{m} (P_{M} \psi_{x}), P_{M} \psi_{x} \rangle dt
\end{equation}

\begin{equation}\label{6.21}
 \lesssim \| \partial_{x} A_{x} \|_{L_{t,x}^{2}} \| P_{\frac{M}{4} \leq \cdot \leq 4M} \psi_{x} \|_{L_{t,x}^{4}} \| P_{M} \psi_{x} \|_{L_{t,x}^{4}} \lesssim \| \partial_{x} A_{x} \|_{L_{t,x}^{2}} \alpha(M)^{2}.
\end{equation}

\noindent Integrating by parts

\begin{equation}\label{6.22}
 - \int_{I} \int (P_{\leq \frac{M}{100}} A_{m}) \partial_{m} |P_{M} \psi_{x}|^{2}dx dt \lesssim \| \partial_{x} A_{x} \|_{L_{t,x}^{2}} \alpha(M)^{2}.
\end{equation}

\noindent Therefore,

\begin{equation}\label{6.23}
\aligned
 \| P_{M} \psi_{x}(t) \|_{L^{2}(\mathbf{R}^{2})}^{2} \leq \beta(M)^{2} + C(E_{0}) \alpha(M)^{2} \| \nabla A_{x} \|_{L_{t,x}^{2}} \\+ C(E_{0}) \| (P_{M} \bar{\psi}_{x}) \psi_{x} \|_{L_{t,x}^{2}} (\| \partial_{x} A_{x} \|_{L_{t,x}^{2}} + \| A_{t} \|_{L_{t,x}^{2}}  + \| A_{x} \|_{L_{t,x}^{4}}^{2} + \| \psi_{x} \|_{L_{t,x}^{4}}^{2}).
\endaligned
\end{equation}

\noindent This proves the theorem assuming

\begin{equation}\label{6.23.1}
 \| \partial_{x} A_{x}(0) \|_{L_{t,x}^{2}} + \| A_{t}(0) \|_{L_{t,x}^{2}} + \| A_{x}(0) \|_{L_{t,x}^{4}}^{2} + \| \psi_{x}(0) \|_{L_{t,x}^{4}}^{2} \lesssim \epsilon^{2}.
\end{equation}

\begin{equation}\label{6.24}
 \| \psi_{x} \|_{L_{t,x}^{4}}^{2} \lesssim \sum_{k} \alpha(k)^{2} \lesssim \epsilon^{2}.
\end{equation}

\noindent Combining $(\ref{6.5})$ and $(\ref{6.7})$,

\begin{equation}\label{6.25}
\| A_{x}(s) \|_{L_{s}^{2} L_{x}^{\infty}([0, \infty) \times \mathbf{R}^{2})} \lesssim_{E_{0}} 1.
\end{equation}

\noindent \textbf{Remark:} For the rest of this section $A \lesssim B$ means $A \lesssim_{E_{0}} B$.

\begin{lemma}\label{l6.5}
 For any $k \geq 0$,

\begin{equation}\label{6.26}
 \| \partial_{x}^{k} \psi_{x}(s) \|_{L_{x}^{4}(\mathbf{R}^{2})} \lesssim_{k} s^{-k/2} \| \psi_{x}(0) \|_{L_{x}^{4}(\mathbf{R}^{2})}.
\end{equation}

\begin{equation}\label{6.27}
 \| \partial_{x} \psi_{x}(s) \|_{L_{s}^{2} L_{t, x}^{4}}^{2} \lesssim \sum_{k} \alpha(k)^{2} \lesssim \| \psi_{x}(0) \|_{L_{t,x}^{4}}^{2}.
\end{equation}
\end{lemma}

\noindent \emph{Proof:} This is proved by Duhamel's principle.

\begin{equation}\label{6.28}
 \psi_{x}(s) = e^{s \Delta} \psi_{x}(0) + \int_{0}^{s} e^{(s - s') \Delta} [2i \partial_{l} (A_{l} \psi_{x}) - (A_{x} \cdot A_{x} + i \partial_{l} A_{l}) \psi_{x} + i Im(\psi_{x} \bar{\psi}_{l}) \psi_{l}](s') ds'.
\end{equation}

\begin{equation}\label{6.29}
 \| \psi_{x}(s) \|_{L_{s}^{\infty} L_{x}^{4}} \lesssim \| \psi_{x}(0) \|_{L_{x}^{4}} + \| \psi_{x} \|_{L_{s}^{2} L_{x}^{\infty}}^{2} \| \psi_{x} \|_{L_{s}^{\infty} L_{x}^{4}} + \| A_{x} \|_{L_{s}^{2} L_{x}^{\infty}}^{2} \| \psi_{x} \|_{L_{s}^{\infty} L_{x}^{4}}
\end{equation}

\begin{equation}\label{6.30}
 + \| \partial_{x} A \|_{L_{s}^{1} L_{x}^{\infty}} \| \psi_{x} \|_{L_{s}^{\infty} L_{x}^{4}} + C(\delta, E_{0}) \| A_{x} \|_{L_{s}^{2} L_{x}^{\infty}} \| \psi_{x} \|_{L_{s}^{\infty} L_{x}^{4}} + \delta \| \psi_{x} \|_{L_{s}^{\infty} L_{x}^{4}}.
\end{equation}

\noindent The last inequality follows from corollary $\ref{c6.2}$ and splitting

\begin{equation}\label{6.31}
 \int_{0}^{s} e^{(s - s') \Delta} \partial_{l} (A_{l} \psi_{x})(s') ds' = \int_{0}^{(1 - \delta) s} e^{(s - s') \Delta} \partial_{l} (A_{l} \psi_{x})(s') ds' + \int_{(1 - \delta)s}^{s} e^{(s - s') \Delta} \partial_{l} (A_{l} \psi_{x})(s') ds'.
\end{equation}

\begin{equation}\label{6.32}
 \int_{(1 - \delta)s}^{s} \frac{1}{(s - s')^{1/2}} \frac{1}{(s')^{1/2}} ds' \lesssim \delta^{1/2}.
\end{equation}

\begin{equation}\label{6.33}
 \int_{0}^{(1 - \delta)s} \frac{1}{(s - s')^{1/2}} f(s') ds' \lesssim C(\delta) \| f \|_{L_{s}^{2}}.
\end{equation}

\noindent By theorem $\ref{t6.1}$ and corollary $\ref{c6.2}$, after partitioning $[0, \infty)$ into finitely many pieces and iterating, $(\ref{6.29})$ and $(\ref{6.30})$ imply

\begin{equation}\label{6.34}
\aligned
 \| \psi_{x}(s) \|_{L_{x}^{4}(\mathbf{R}^{2})} &\lesssim \| \psi_{x}(0) \|_{L_{x}^{4}(\mathbf{R}^{2})}, \\
\| \psi_{x}(s) \|_{L_{t, x}^{4}(I \times \mathbf{R}^{2})} &\lesssim \| \psi_{x}(0) \|_{L_{t, x}^{4}(I \times \mathbf{R}^{2})}.
\endaligned
\end{equation}

\noindent Likewise since the kernel of $\partial_{x}^{k} e^{(s - s') \Delta}$ has $L^{1}$ norm bounded by $\frac{1}{(s - s')^{k/2}}$,

\begin{equation}\label{6.35}
\| \int_{0}^{(1 - \delta) s} e^{(s - s') \Delta} \partial_{x}^{k} [2i \partial_{l} (A_{l} \psi_{x}) - (A_{x} \cdot A_{x} + i \partial_{l} A_{l}) \psi_{x} + i Im(\psi_{x} \bar{\psi}_{l}) \psi_{l}](s') ds' \|_{L_{x}^{4}(\mathbf{R}^{2})} \lesssim_{k} s^{-k/2} \| \psi_{x}(0) \|_{L_{x}^{4}}.
\end{equation}

\noindent Next, theorem $\ref{t6.1}$, corollary $\ref{c6.2}$, and an induction imply

\begin{equation}\label{6.36}
\| \int_{(1 - \delta)s}^{s} \partial_{x}^{k} e^{(s - s') \Delta}[- (A_{x} \cdot A_{x} + i \partial_{l} A_{l}) \psi_{x} + i Im(\psi_{x} \bar{\psi}_{l}) \psi_{l}](s') ds' \|_{L_{x}^{4}}
\end{equation}

\begin{equation}\label{6.37}
 \lesssim \int_{(1 - \delta)s}^{s} \frac{1}{(s - s')^{1/2}} \| \partial_{x}^{k - 1} [- (A_{x} \cdot A_{x} + i \partial_{l} A_{l}) \psi_{x} + i Im(\psi_{x} \bar{\psi}_{l}) \psi_{l}](s') \|_{L_{x}^{4}} ds' \lesssim_{k} s^{-k/2} \| \psi_{x}(0) \|_{L_{x}^{4}}.
\end{equation}

\noindent Finally,

\begin{equation}\label{6.38}
\| \int_{(1 - \delta)s}^{s} \partial_{x}^{k} e^{(s - s') \Delta} [2i \partial_{l} (A_{l} \psi_{x})](s')ds' \|_{L_{x}^{4}} \lesssim_{E_{0},k} s^{-k/2} \delta + \delta^{1/2} s^{-k/2} \| s^{k/2} \partial_{x}^{k} \psi_{x}(s) \|_{L_{s}^{\infty} L_{x}^{4}}.
\end{equation}

\noindent Combining $(\ref{6.35})$, $(\ref{6.37})$, $(\ref{6.38})$, and $\| \partial_{x}^{k} e^{s \Delta} \psi_{x}(0) \|_{L_{x}^{4}} \lesssim s^{-k/2} \| \psi_{x}(0) \|_{L_{x}^{4}}$ proves $(\ref{6.26})$.\vspace{5mm}

\noindent Now to prove $(\ref{6.27})$. Estimate

\begin{equation}\label{6.39}
 \| \int_{s/2}^{s} e^{(s - s') \Delta} \partial_{x} [2i \partial_{l} (A_{l} \psi_{x}) - (A_{x} \cdot A_{x} + i \partial_{l} A_{l}) \psi_{x} + i Im(\psi_{x} \bar{\psi}_{l}) \psi_{l}](s') ds' \|_{L_{s}^{2} L_{x}^{4}([2^{j}, 2^{j + 1}] \times \mathbf{R}^{2})}
\end{equation}

\noindent By theorem $\ref{t6.1}$, corollary $\ref{c6.2}$, and $(\ref{6.26})$,

\begin{equation}\label{6.40}
 \lesssim (\| \nabla A \|_{L_{s}^{1} L_{x}^{\infty}([2^{j - 2}, 2^{j}] \times \mathbf{R}^{2})}^{1/2} + \| \psi_{x} \|_{L_{s}^{2} L_{x}^{\infty}([2^{j - 2}, 2^{j}] \times \mathbf{R}^{2})} + \| A_{x} \|_{L_{s}^{2} L_{x}^{\infty}([2^{j - 2}, 2^{j}] \times \mathbf{R}^{2})}) \| \psi_{x} \|_{L_{s}^{\infty} L_{x}^{4}}.
\end{equation}

\noindent Next, by Sobolev embedding, theorem $\ref{t6.1}$, and corollary $\ref{c6.2}$, and Holder's inequality,

\begin{equation}\label{6.41}
 \| \int_{0}^{s/2} e^{(s - s') \Delta} \partial_{x} [2i \partial_{l} (A_{l} \psi_{x}) - (A_{x} \cdot A_{x} + i \partial_{l} A_{l}) \psi_{x} + i Im(\psi_{x} \bar{\psi}_{l}) \psi_{l}](s') ds' \|_{L_{s}^{2} L_{x}^{4}([2^{j}, 2^{j + 1}] \times \mathbf{R}^{2})}
\end{equation}

\begin{equation}\label{6.42}
 \lesssim \| \psi_{x} \|_{L_{s}^{\infty} L_{x}^{4}} 2^{-j/2} \sum_{k \leq j} 2^{k/2} (\| A_{x} \|_{L_{s,x}^{4}([2^{k}, 2^{k + 1}] \times \mathbf{R}^{2})}^{2} + \| \partial_{x} A \|_{L_{s,x}^{2}([2^{k}, 2^{k + 1}] \times \mathbf{R}^{2})} + \| \psi_{x} \|_{L_{s,x}^{4}([2^{k}, 2^{k + 1}] \times \mathbf{R}^{2})}^{2}).
\end{equation}

\noindent Combining $(\ref{6.40})$ and $(\ref{6.42})$ implies

\begin{equation}\label{6.43}
 \| \int_{0}^{s} e^{(s - s') \Delta} \partial_{k} [2i \partial_{l} (A_{l} \psi_{x}) - (A_{x} \cdot A_{x} + i \partial_{l} A_{l}) \psi_{x} + i Im(\psi_{x} \bar{\psi}_{l}) \psi_{l}](s') ds' \|_{L_{s}^{2} L_{x}^{4}} \lesssim_{E_{0}} \| \psi_{x}(0) \|_{L_{x}^{4}}.
\end{equation}

\begin{equation}\label{6.44}
 \| \nabla e^{s \Delta} \psi_{x}(0) \|_{L_{t,x}^{4}}^{2} \lesssim \sum_{2^{2k} \leq \frac{1}{s}} 2^{k} \| P_{k} \psi_{x}(0) \|_{L_{t,x}^{4}}^{2} + \sum_{2^{2k} > \frac{1}{s}} \frac{2^{k}}{(s 2^{2k})^{3}} \| P_{k} \psi_{x}(0) \|_{L_{t,x}^{4}}^{2}.
\end{equation}

\begin{equation}\label{6.45}
 \| \nabla e^{s \Delta} \psi_{x}(0) \|_{L_{s}^{2} L_{t,x}^{4}([2^{-2j}, 2^{-2j + 2}] \times I \times \mathbf{R}^{2}}^{2} \lesssim \sum_{k \geq j} 2^{k - j} \| P_{k} \psi_{x}(0) \|_{L_{t,x}^{4}}^{2} + \sum_{k > j} 2^{5(j - k)} \| P_{k} \psi_{x}(0) \|_{L_{t,x}^{4}}^{2}.
\end{equation}

\noindent Therefore,

\begin{equation}\label{6.46}
 \sum_{j} \| \nabla e^{s \Delta} \psi_{x}(0) \|_{L_{s}^{2} L_{t,x}^{4}([2^{-2j}, 2^{-2j + 2}] \times I \times \mathbf{R}^{2}}^{2} \lesssim \sum_{k} \beta(k)^{2} \lesssim \| \psi_{x}(0) \|_{L_{t,x}^{4}}^{2}.
\end{equation}

\noindent This gives $(\ref{6.27})$. $\Box$

\begin{corollary}\label{c6.6}
\begin{equation}\label{6.47}
\| A_{x}(s) \|_{L_{t,x}^{4}} \lesssim \epsilon.
\end{equation}

\end{corollary}

\noindent \emph{Proof:} This follows from lemma $\ref{l6.5}$, theorem $\ref{t6.1}$, corollary $\ref{c6.2}$, and the formula

\begin{equation}\label{6.48}
 A_{x}(s) = -\int_{s}^{\infty} Im(\bar{\psi}_{x} (\partial_{l} \psi_{l} + i A_{l} \psi_{l}))(r) dr.
\end{equation}

\begin{equation}\label{6.48.1}
\| A_{x} \|_{L_{t,x}^{4}} \lesssim \| \psi_{x} \|_{L_{s}^{2} L_{x}^{\infty}} (\| \partial_{x} \psi_{x} \|_{L_{s}^{2} L_{t,x}^{4}} + \| A_{x} \|_{L_{s}^{2} L_{x}^{\infty}} \| \psi_{x} \|_{L_{s}^{\infty} L_{x}^{4}}) \lesssim \| \psi_{x}(0) \|_{L_{t,x}^{4}}.
\end{equation}

\noindent $\Box$

\begin{theorem}\label{t6.7}
 \begin{equation}\label{6.49}
  \| P_{k} A_{x}(s) \|_{L_{x}^{1}(\mathbf{R}^{2})} \lesssim 2^{-k}.
 \end{equation}
\end{theorem}

\noindent \emph{Proof:} By Bernstein's inequality

\begin{equation}\label{6.50}
 2^{2k} \| P_{k} A_{x}(s) \|_{L_{x}^{1}(\mathbf{R}^{2})} \lesssim \int_{s}^{\infty} \| \nabla^{2} \psi_{x}(r) \|_{L_{x}^{2}} \| \nabla \psi_{x}(r) \|_{L_{x}^{2}} + \| \psi_{x}(r) \|_{L_{x}^{2}} \| \nabla^{3} \psi_{x}(r) \|_{L_{x}^{2}} dr
\end{equation}

\begin{equation}\label{6.51}
 + \int_{s}^{\infty} \| \nabla^{2} \psi_{x} \|_{L^{2}} \| A_{x} \|_{L^{\infty}} \| \psi_{x} \|_{L^{2}} + \| \psi_{x} \|_{L^{2}} \| \nabla^{2} A_{x} \|_{L_{x}^{2}} \| \psi_{x} \|_{L_{x}^{\infty}} + \| \psi_{x} \|_{L_{x}^{2}} \| A_{x} \|_{L_{x}^{2}} \| \nabla^{2} \psi_{x} \|_{L_{x}^{\infty}} dr
\end{equation}

\begin{equation}\label{6.52}
 \lesssim \int_{s}^{\infty} r^{-3/2} dr \lesssim s^{-1/2}.
\end{equation}

\noindent The first inequality in $(\ref{6.52})$ follows from theorem $\ref{t6.1}$ and corollary $\ref{c6.2}$. So for $s > 2^{-2k}$,

\begin{equation}\label{6.53}
 \| P_{k} A_{x}(s) \|_{L_{x}^{1}(\mathbf{R}^{2})} \lesssim 2^{-k}.
\end{equation}

\noindent For $s < 2^{-2k}$ Holder's inequality and $(\ref{6.53})$ imply

\begin{equation}\label{6.54}
 \| P_{k} A_{x}(s) \|_{L_{x}^{1}(\mathbf{R}^{2})} \lesssim \int_{s}^{2^{-2k}} \| |\psi_{x}| |\partial_{l} \psi_{l} + i A_{l} \psi_{l}| \|_{L_{x}^{1}} dr
\end{equation}

\begin{equation}\label{6.55}
 \lesssim 2^{-k} \| \psi_{x} \|_{L_{s}^{\infty} L_{x}^{2}} \| \partial_{x} \psi_{x} \|_{L_{s,x}^{2}} + 2^{-k} \| \psi_{x} \|_{L_{s}^{\infty} L_{x}^{2}}^{2} \| A_{x} \|_{L_{s}^{2} L_{x}^{\infty}} + 2^{-k} \lesssim 2^{-k}.
\end{equation}

\noindent $\Box$

\begin{lemma}\label{l6.8}
 \begin{equation}\label{6.56}
  \| \psi_{t} \|_{L_{s}^{2} L_{t,x}^{4}}^{2} \lesssim \sum_{k} \beta(k)^{2} \lesssim \epsilon^{2}.
 \end{equation}

\end{lemma}

\noindent \emph{Proof:}

\begin{equation}\label{6.57}
 \psi_{t}(0) = i \partial_{l} \psi_{l}(0) - A_{l}(0) \psi_{l}(0).
\end{equation}

\noindent As in lemma $\ref{l6.5}$

\begin{equation}\label{6.58}
 \| e^{s \Delta} (\partial_{l} \psi_{l}(0)) \|_{L_{s}^{2} L_{t,x}^{4}} \lesssim \epsilon.
\end{equation}

\noindent By Sobolev embedding and theorem $\ref{t6.7}$

\begin{equation}\label{6.59}
 \| P_{N}(A_{l}(0) \psi_{l}(0)) \|_{L_{t,x}^{4}} \lesssim N \| P_{\frac{N}{20} \leq \cdot \leq 20N} \psi_{l}(0) \|_{L_{t,x}^{4}} \| A_{l}(0) \|_{L_{x}^{2}}
\end{equation}

\begin{equation}\label{6.60}
 + N^{1/2} \| P_{\leq \frac{N}{20}} \psi_{l}(0) \|_{L_{t}^{4} L_{x}^{\infty}} \| P_{\frac{N}{4} \leq \cdot \leq 4N} A_{l}(0) \|_{L_{x}^{2}}
\end{equation}

\begin{equation}\label{6.61}
 + N^{3/2} \sum_{k > 0} \| P_{2^{k} N} \psi_{x}(0) \|_{L_{t,x}^{4}} \| P_{2^{k} N} A_{x}(0) \|_{L_{x}^{4/3}}
\end{equation}

\begin{equation}\label{6.62}
 \lesssim N \alpha(N) + N \sum_{k \leq 0} 2^{k/2} \alpha(2^{k} N) + N \sum_{k > 0} 2^{-k/2} \alpha(2^{k} N) \lesssim N \alpha(N).
\end{equation}

\noindent This implies

\begin{equation}\label{6.63}
 \| e^{s \Delta} \psi_{t}(0) \|_{L_{s}^{2} L_{t,x}^{4}} + \| s^{1/2} e^{s \Delta} \psi_{t}(0) \|_{L_{t,x}^{4}} \lesssim \epsilon.
\end{equation}

\begin{equation}\label{6.64}
 \psi_{t}(s) = e^{s \Delta} \psi_{t}(0) + \int_{0}^{s} e^{(s - s') \Delta} [2i \partial_{l} (A_{l} \psi_{t}) - (A_{x} \cdot A_{x} + i \partial_{l} A_{l}) \psi_{t} + i Im(\psi_{t} \bar{\psi}_{l}) \psi_{l}](s') ds',
\end{equation}

\noindent Making an argument identical to the proof of lemma $\ref{l6.5}$ proves the lemma. $\Box$

\begin{corollary}\label{c6.9}
\begin{equation}\label{6.65}
\| A_{t}(s) \|_{L_{t,x}^{2}} \lesssim \epsilon^{2}.
\end{equation}

\end{corollary}

\noindent \emph{Proof:}

\begin{equation}\label{6.66}
 \| A_{t} \|_{L_{t,x}^{2}} \lesssim \| \psi_{t} \|_{L_{s}^{2} L_{t,x}^{4}} (\| \partial_{x} \psi_{x} \|_{L_{s}^{2} L_{t,x}^{4}} + \| A_{x} \|_{L_{s}^{2} L_{x}^{\infty}} \| \psi_{x} \|_{L_{s}^{\infty} L_{t,x}^{4}}) \lesssim \epsilon^{2}.
\end{equation}

\noindent Recall the choice of frequency envelope

\begin{equation}\label{6.66.1}
 \alpha(k) = \sum_{j} 2^{-\delta |j - k|} \| P_{j} \psi_{x}(0) \|_{L_{t,x}^{4}}.
\end{equation}

\noindent Let

\begin{equation}\label{6.66.2}
\alpha(t,k) = \sum_{j} 2^{-\delta|j - k|} \| P_{j} \psi_{x}(t,0) \|_{L_{x}^{4}}.
\end{equation}

\begin{equation}\label{6.66.3}
(\int \alpha(t,k)^{4} dt)^{1/4} \lesssim \sum_{j} 2^{-\delta|j - k|} \| P_{j} \psi_{x}(0) \|_{L_{t,x}^{4}} = \alpha(k).
\end{equation}

\begin{equation}\label{6.66.4}
(\int (\sum_{k} \| P_{k} \psi_{x}(t,0) \|_{L_{x}^{4}}^{2})^{2} dt)^{1/2} \lesssim \sum_{k} (\int \| P_{k} \psi_{x}(t,0) \|_{L_{x}^{4}}^{4} dt)^{1/2} \lesssim \sum \alpha(k)^{2} \lesssim \epsilon^{2}.
\end{equation}

\begin{theorem}\label{t6.10}
\begin{equation}\label{6.67}
 \| P_{k} \psi_{x}(s) \|_{L_{x}^{4}} \lesssim (1 + s 2^{2k})^{-4} \alpha(t,k).
\end{equation}

\end{theorem}

\noindent \emph{Proof:} We start by proving $\| P_{k} \psi_{x}(s) \|_{L_{t,x}^{4}} \lesssim \alpha(k)$.

\begin{equation}\label{6.69}
 \| e^{s \Delta} P_{k} \psi_{x}(0) \|_{L_{x}^{4}} \lesssim (1 + s 2^{2k})^{-4} \alpha(t,k).
\end{equation}

\noindent Make the bootstrap assumption

\begin{equation}\label{6.68}
 \| P_{k} \psi_{x}(s) \|_{L_{t,x}^{4}} \leq C \alpha(t,k).
\end{equation}

\begin{equation}\label{6.75}
 \| \int_{0}^{(1 - \delta)s} e^{(s - s') \Delta} P_{k}[2i \partial_{l}(A_{l} \psi_{x}) - (A_{l} A_{l} + i \partial_{l} A_{l}) \psi_{x} + i Im(\psi_{x} \bar{\psi}_{l}) \psi_{l}] ds' \|_{L_{x}^{4}}
\end{equation}

\begin{equation}\label{6.76}
\aligned
 \lesssim   e^{-\delta s 2^{2k}} \| P_{k - 5 \leq \cdot \leq k + 5} \psi_{x} \|_{L_{s}^{\infty} L_{x}^{4}} [s^{1/2} 2^{k} \| A_{x} \|_{L_{s}^{2} L_{x}^{\infty}} \\ +  \| A_{x} \|_{L_{s}^{2} L_{x}^{\infty}}^{2} + \| \partial_{x} A_{x} \|_{L_{s}^{1} L_{x}^{\infty}} + \| \psi_{x} \|_{L_{s}^{2} L_{x}^{\infty}}^{2}]
\endaligned
\end{equation}

\begin{equation}\label{6.77}
\aligned
+ e^{-\delta s 2^{2k}} (\sum_{j \leq k - 5} 2^{j/2} \| P_{j} \psi_{x} \|_{L_{s}^{\infty} L_{x}^{4}}) [s^{1/2} 2^{k} \| P_{\geq k - 5} A_{x} \|_{L_{s}^{2} L_{x}^{4}}  \\
 + \| P_{\geq k - 5} A_{x} \|_{L_{s}^{2} L_{x}^{4}} \| A_{x} \|_{L_{s}^{2} L_{x}^{\infty}} + \| \partial_{x} P_{\geq k - 5} A_{x} \|_{L_{s}^{1} L_{x}^{\infty}} + \| P_{\geq k - 5} \psi_{x} \|_{L_{s}^{2} L_{x}^{4}} \| \psi_{x} \|_{L_{s}^{2} L_{x}^{\infty}}]
\endaligned
\end{equation}

\begin{equation}\label{6.78}
\aligned
 + e^{-s \delta 2^{2k}} s^{1/2} 2^{k} \sum_{j \geq k} \| P_{j} \psi_{x} \|_{L_{s}^{\infty} L_{t,x}^{4}} [\| P_{j} A_{x} \|_{L_{s,x}^{2}}  \\
+ \| P_{\geq j} A_{x} \|_{L_{s,x}^{2}} \| A_{x} \|_{L_{s}^{2} L_{x}^{\infty}} + \| \partial_{x} P_{\geq j} A_{x} \|_{L_{s}^{1} L_{x}^{2}} + \|  P_{\geq j} \psi_{x} \|_{L_{s,x}^{2}} \| \psi_{x} \|_{L_{s}^{2} L_{x}^{\infty}}].
\endaligned
\end{equation}

\noindent Next, estimate

\begin{equation}\label{6.70}
\| \int_{(1 - \delta)s}^{s} e^{(s - s') \Delta} P_{k} [2i \partial_{l}(A_{l} \psi_{x}) - (A_{l} A_{l} + i \partial_{l} A_{l}) \psi_{x} + i Im(\psi_{x} \bar{\psi}_{l}) \psi_{l}] ds' \|_{L_{x}^{4}},
\end{equation}

\noindent By Sobolev embedding and integration,

\begin{equation}\label{6.71}
\aligned
 \lesssim   \| P_{k - 5 \leq \cdot \leq k + 5} \psi_{x} \|_{L_{s}^{\infty} L_{x}^{4}} [(\int_{0}^{\delta s} e^{-s' 2^{2k}} 2^{2k} ds')^{1/2} \delta^{1/2} \| s^{1/2} A_{x} \|_{L_{s,x}^{\infty}} \\ + \delta \| s^{1/2} A_{x} \|_{L_{s,x}^{\infty}}^{2} + \delta \| s \partial_{x} A_{x} \|_{L_{s,x}^{\infty}} + \delta \| s^{1/2} \psi_{x} \|_{L_{s,x}^{\infty}}^{2}
\endaligned
\end{equation}

\begin{equation}\label{6.72}
\aligned
+ (\sum_{j \leq k - 5} 2^{j/2} \| P_{j} \psi_{x} \|_{L_{s}^{\infty} L_{x}^{4}}) [\delta^{1/2} \| s^{1/2} P_{\geq k - 5} A_{x} \|_{L_{s}^{\infty} L_{x}^{4}} (\int_{0}^{\delta s} e^{-s' 2^{2k}} 2^{2k} ds')^{1/2} \\
 + \delta \| s^{1/2} P_{\geq k - 5} A_{x} \|_{L_{s}^{\infty} L_{x}^{4}} \| s^{1/2} A_{x} \|_{L_{s,x}^{\infty}} + \delta \| s \partial_{x} P_{\geq k - 5} A_{x} \|_{L_{s}^{\infty} L_{x}^{4}} + \delta \| s^{1/2} P_{\geq k - 5} \psi_{x} \|_{L_{s}^{\infty} L_{x}^{4}} \| s^{1/2} \psi_{x} \|_{L_{s,x}^{\infty}}]
\endaligned
\end{equation}

\begin{equation}\label{6.73}
\aligned
 + 2^{k} \sum_{j \geq k} \| P_{j} \psi_{x} \|_{L_{s}^{\infty} L_{x}^{4}} [\delta^{1/2} \| s^{1/2} P_{j} A_{x} \|_{L_{s}^{\infty} L_{x}^{2}} (\int_{0}^{\delta s} e^{-s' 2^{2k}} 2^{2k} ds')^{1/2} \\
+ \delta \| s^{1/2} P_{\geq j} A_{x} \|_{L_{s}^{\infty} L_{x}^{2}} \| s^{1/2} A_{x} \|_{L_{s,x}^{\infty}} + \delta \| s \partial_{x} P_{\geq j} A_{x} \|_{L_{s}^{\infty} L_{x}^{2}} + \delta \| s^{1/2} P_{\geq j} \psi_{x} \|_{L_{s}^{\infty} L_{x}^{2}} \| s^{1/2} \psi_{x} \|_{L_{s,x}^{\infty}}
\endaligned
\end{equation}

\begin{equation}\label{6.74}
 \lesssim C \alpha(t, k) \delta^{1/2}.
\end{equation}

\noindent The last inequality follows from Bernstein's inequality, the bootstrap assumption, theorem $\ref{t6.1}$, and corollary $\ref{c6.2}$. Partitioning $[0, \infty)$ into finitely many intervals $I_{j}$ for each $t$ such that

\begin{equation}\label{6.74.1}
\| A_{x} \|_{L_{s}^{2} L_{x}^{\infty}(I_{j} \times \mathbf{R}^{2})} + \| \partial_{x} A_{x} \|_{L_{s}^{1} L_{x}^{\infty}(I_{j} \times \mathbf{R}^{2})} + \| \psi_{x} \|_{L_{s}^{2} L_{x}^{\infty}(I_{j} \times \mathbf{R}^{2})} + \| \partial_{x} A_{x} \|_{L_{s,x}^{2}(I_{j} \times \mathbf{R}^{2})}
\end{equation}

\noindent is small on each $I_{j}$ and iterating,

\begin{equation}\label{6.75}
 \| P_{k} \psi_{x}(s,t) \|_{L_{x}^{4}} \lesssim \alpha(t,k).
\end{equation}

\noindent This in turn implies

\begin{equation}\label{6.81}
 \| P_{k} \psi_{x}(s) \|_{L_{t,x}^{4}} \lesssim \alpha(k).
\end{equation}

\noindent To prove $(\ref{6.67})$ it only remains to consider $s > 2^{-2k}$. $e^{-s \delta 2^{2k}} s^{1/2} 2^{k} \lesssim_{\delta} (1 + s 2^{2k})^{-4}$, which takes care of $(\ref{6.76})$, $(\ref{6.77})$, and $(\ref{6.78})$. Now make the bootstrap assumption

\begin{equation}\label{6.82}
 \| P_{k} \psi_{x}(s) \|_{L_{t,x}^{4}} \leq C(1 + s 2^{2k})^{-4} \alpha(k).
\end{equation}

\noindent Plugging this in to $(\ref{6.71})$ and $(\ref{6.73})$

\begin{equation}\label{6.83}
 (\ref{6.71}) + (\ref{6.73}) \lesssim C \delta^{1/2} (1 + s 2^{2k})^{-4} \alpha(k).
\end{equation}

\noindent By theorem $\ref{t6.1}$, corollary $\ref{c6.2}$, and Bernstein's inequality,

\begin{equation}\label{6.83.1}
\delta^{1/2} \| s^{1/2} P_{\geq k - 5} A_{x} \|_{L_{s}^{\infty} L_{x}^{4}} + \delta \| s^{1/2} P_{\geq k - 5} A_{x} \|_{L_{s}^{\infty} L_{x}^{4}} \| s^{1/2} A_{x} \|_{L_{s,x}^{\infty}}
\end{equation}

\begin{equation}\label{6.83.2}
+ \delta \| s \partial_{x} P_{\geq k - 5} A_{x} \|_{L_{s,x}^{\infty}} + \delta \| s^{1/2} P_{\geq k - 5} \psi_{x} \|_{L_{s}^{\infty} L_{x}^{4}} \| s^{1/2} \psi_{x} \|_{L_{s,x}^{\infty}}
\end{equation}

\begin{equation}\label{6.83.3}
 \lesssim s^{-4} 2^{-17k/2} \delta^{1/2} \| s^{9/2} \partial_{x}^{9} A_{x} \|_{L_{s}^{\infty} L_{x}^{2}}^{1/2} \| s^{9/2} \partial_{x}^{8} A_{x} \|_{L_{s,x}^{\infty}}^{1/2} 
\end{equation}

\begin{equation}\label{6.83.5}
+ s^{-4} 2^{-17k/2} \delta^{1/2} \| s^{9/2} \partial_{x}^{9} A_{x} \|_{L_{s}^{\infty} L_{x}^{2}}^{1/2} \| s^{9/2} \partial_{x}^{8} A_{x} \|_{L_{s,x}^{\infty}}^{1/2} \| s^{1/2} A_{x} \|_{L_{s,x}^{\infty}}
\end{equation}

\begin{equation}\label{6.83.4}
+ s^{-4} 2^{-17k/2} \delta^{1/2} \| s^{5} \partial_{x}^{10} A_{x} \|_{L_{s}^{\infty} L_{x}^{2}}^{1/2} \| s^{5} \partial_{x}^{10} A_{x} \|_{L_{s,x}^{\infty}}^{1/2} 
\end{equation}

\begin{equation}\label{6.83.6}
 + s^{-4} 2^{-17k/2} \delta^{1/2} \| s^{9/2} \partial_{x}^{9} \psi_{x} \|_{L_{s}^{\infty} L_{x}^{2}}^{1/2} \| s^{9/2} \partial_{x}^{8} \psi_{x} \|_{L_{s,x}^{\infty}}^{1/2} \| s^{1/2} \psi_{x} \|_{L_{s,x}^{\infty}} \lesssim_{E_{0}} \delta^{1/2} s^{-4} 2^{-17k/2}.
\end{equation}

\noindent Since $\| P_{k} \psi_{x}(s) \|_{L_{t,x}^{4}} \lesssim \alpha(k)$, when $s > 2^{-2k}$

\begin{equation}\label{6.84}
 (\ref{6.72}) \lesssim (1 + s 2^{2k})^{-4}.
\end{equation}

\noindent This completes the proof of the theorem. $\Box$

\begin{corollary}\label{c6.11}
\begin{equation}\label{6.85}
 2^{k} \| P_{k} A_{x}(s) \|_{L_{t,x}^{2}} \lesssim \epsilon \alpha(k).
\end{equation}
\end{corollary}

\noindent \emph{Proof:}

\begin{equation}\label{6.86}
 A_{x}(s) = -\int_{s}^{\infty} Im(\bar{\psi}_{x} (\partial_{l} \psi_{l} + i A_{l} \psi_{l}))(r) dr.
\end{equation}

\begin{equation}\label{6.87}
 \| P_{k} A_{x}(s) \|_{L_{t,x}^{2}} \lesssim \int_{s}^{\infty} \| P_{k - 5 \leq \cdot \leq k + 5} \psi_{x}(r) \|_{L_{t,x}^{4}} (\| \partial_{x} \psi_{x}(r) \|_{L_{t,x}^{4}} + \| A_{x}(r) \|_{L_{x}^{\infty}} \| \psi_{x}(r) \|_{L_{t,x}^{4}}) dr
\end{equation}

\begin{equation}\label{6.88}
 + \int_{s}^{\infty} \| P_{\leq k - 5}  \psi_{x}(r) \|_{L_{t,x}^{4}} \| \partial_{x} P_{\geq k - 5} \psi_{x}(r) \|_{L_{t,x}^{4}} dr
\end{equation}

\begin{equation}\label{6.89}
 + \int_{s}^{\infty} \| P_{\leq k - 5} \psi_{x}(r) \|_{L_{t,x}^{4}} \| P_{\geq k - 5} \psi_{x}(r) \|_{L_{t,x}^{4}} \| A_{x}(r) \|_{L_{x}^{\infty}} + \| P_{\leq k - 5} \psi_{x}(r) \|_{L_{t}^{4} L_{x}^{\infty}}^{2} \| P_{\geq k - 5} A_{x}(r) \|_{L_{x}^{2}}
\end{equation}

\begin{equation}\label{6.90}
\aligned
 + \sum_{j \geq k + 5}  \int_{s}^{\infty} 2^{j} \| P_{j} \psi_{x}(r) \|_{L_{t,x}^{4}}^{2} + \| P_{j} \psi_{x}(r) \|_{L_{t,x}^{4}} \| P_{\geq j} \psi_{x}(r) \|_{L_{t,x}^{4}} \| A_{x}(r) \|_{L_{x}^{\infty}} \\ + \| P_{j} \psi_{x}(r) \|_{L_{t,x}^{4}} \| \psi_{x}(r) \|_{L_{t,x}^{4}} \| P_{\geq j} A_{x}(r) \|_{L_{x}^{\infty}} dr.
\endaligned
\end{equation}

\begin{equation}\label{6.91}
(\ref{6.87}) \lesssim (\int_{0}^{\infty} (1 + s 2^{2k})^{-8} \alpha(k) ds)^{1/2} (\| \partial_{x} \psi_{x} \|_{L_{s}^{2} L_{t,x}^{4}} + \| A_{x} \|_{L_{s}^{2} L_{x}^{\infty}} \| \psi_{x} \|_{L_{s}^{\infty} L_{t,x}^{4}}) \lesssim 2^{-k} \alpha(k) \epsilon.
\end{equation}

\begin{equation}\label{6.92}
 (\ref{6.88}) \lesssim \sum_{j \geq k - 5} \epsilon \int_{0}^{\infty} (1 + s 2^{2j})^{-4} \alpha(j) ds \lesssim 2^{-k} \epsilon \alpha(k).
\end{equation}

\noindent By Bernstein's inequality, theorem $\ref{t6.1}$, corollary $\ref{c6.2}$,

\begin{equation}\label{6.93}
(\ref{6.89}) \lesssim \epsilon (\sum_{j} \int_{0}^{\infty} (1 + s 2^{2j})^{-8} ds)^{1/2} \| A_{x} \|_{L_{s}^{2} L_{x}^{\infty}} + \alpha(k)^{2} 2^{-k} \| \partial_{x}^{2} A_{x} \|_{L_{s}^{1} L_{x}^{2}} \lesssim 2^{-k} \epsilon \alpha(k).
\end{equation}

\begin{equation}\label{6.94}
(\ref{6.90}) \lesssim \sum_{j \geq k + 5} 2^{j} \epsilon (\int_{0}^{\infty} (1 + s 2^{2j})^{-4} ds) + \epsilon (\int_{0}^{\infty} (1 + s 2^{2j})^{-8} ds)^{1/2} \| A_{x} \|_{L_{s}^{2} L_{x}^{\infty}} \lesssim 2^{-k} \alpha(k) \epsilon.
\end{equation}

\noindent $\Box$\vspace{5mm}

\noindent In conclusion this proves

\begin{equation}\label{6.95}
 \| A_{x} \|_{L_{t,x}^{4}}^{2} + \| \psi_{x} \|_{L_{t,x}^{4}}^{2} + \| \partial_{x} A_{x} \|_{L_{t,x}^{2}} + \| A_{t} \|_{L_{t,x}^{2}} \lesssim \epsilon^{2}.
\end{equation}

\noindent This completes the proof of theorem $\ref{t6.3}$. $\Box$\vspace{5mm}

\noindent Performing an identical calculation to the one performed in the case of the Coulomb gauge, the error involving terms of the form

\begin{equation}\label{6.96}
 P_{M}((P_{\leq \frac{M}{100}} A_{l}) \partial_{l} \psi_{x})
\end{equation}

\noindent is bounded by $C(E_{0}) \epsilon (\alpha(M) + \beta(M))^{2} (\alpha(N) + \beta(N))^{2}$.

\begin{equation}\label{6.97}
 M \| P_{\geq \frac{M}{100}} A_{x} \|_{L_{t,x}^{2}} + \| \partial_{x} A_{x} \|_{L_{t,x}^{2}} + \| A_{t} \|_{L_{t,x}^{2}} + \| A_{x} \|_{L_{t,x}^{4}}^{2} + \| \psi_{x} \|_{L_{t,x}^{4}}^{2} \lesssim \epsilon^{2}.
\end{equation}

\noindent Therefore in the caloric gauge the error is bounded by

\begin{equation}\label{6.98}
 \| (P_{M} \psi_{x})(P_{N} \bar{\psi}_{x}) \|_{L_{t,x}^{2}} \lesssim (\frac{M}{N})^{1/2} (\alpha(M) + \beta(M))(\alpha(N) + \beta(N))
\end{equation}

\begin{equation}\label{6.99}
 + (\frac{M}{N})^{1/2} C(E_{0}) \epsilon \sum_{k} \| (P_{M} \psi_{x})(P_{2^{k} M} \bar{\psi}_{x}) \|_{L_{t,x}^{2}} (\alpha(N) + \beta(N))
\end{equation}

\begin{equation}\label{6.100}
 + (\frac{M}{N})^{1/2} C(E_{0}) \epsilon \sum_{k} \| (P_{N} \psi_{x})(P_{2^{k} N} \bar{\psi}_{x}) \|_{L_{t,x}^{2}} (\alpha(M) + \beta(M)).
\end{equation}

\begin{theorem}\label{t6.12}
 For $M << N$,

\begin{equation}\label{6.101}
 \| (P_{M} \psi_{x})(P_{N} \bar{\psi}_{x}) \|_{L_{t,x}^{2}} \lesssim (\frac{M}{N})^{1/2} (\alpha(M) + \beta(M))(\alpha(N) + \beta(N)).
\end{equation}

\end{theorem}

\noindent \emph{Proof:} From the Morawetz estimates if $M << N$,

\begin{equation}\label{6.102}
 \| (P_{M} \psi_{x})(P_{N} \bar{\psi}_{x}) \|_{L_{t,x}^{2}} \lesssim (\frac{M}{N})^{1/2} (\alpha(M) + \beta(M))(\alpha(N) + \beta(N))
\end{equation}

\begin{equation}\label{6.103}
 + (\frac{M}{N})^{1/2} C(E_{0}) \epsilon \sum_{k} \| (P_{M} \psi_{x})(P_{2^{k} M} \bar{\psi}_{x}) \|_{L_{t,x}^{2}} (\alpha(N) + \beta(N))
\end{equation}

\begin{equation}\label{6.104}
 + (\frac{M}{N})^{1/2} C(E_{0}) \epsilon \sum_{k} \| (P_{N} \psi_{x})(P_{2^{k} N} \bar{\psi}_{x}) \|_{L_{t,x}^{2}} (\alpha(M) + \beta(M)).
\end{equation}

\noindent Therefore,

\begin{equation}\label{6.105}
 \sum_{M, N, M << N} \| (P_{M} \psi_{x})(P_{N} \bar{\psi}_{x}) \|_{L_{t,x}^{2}} \lesssim \sum_{M} (\alpha(M) + \beta(M))^{2}
\end{equation}

\begin{equation}\label{6.106}
 + C(E_{0}) \epsilon \sum_{M \leq N} (\frac{M}{N}))^{1/2} (\sum_{k} \| (P_{M} \psi_{x}) (P_{2^{k} M} \bar{\psi}_{x}) \|_{L_{t,x}^{2}}) (\alpha(N) + \beta(N))
\end{equation}

\begin{equation}\label{6.107}
 + C(E_{0}) \epsilon \sum_{M \leq N} (\frac{M}{N})^{1/2} (\sum_{k} \| (P_{N} \psi_{x})(P_{2^{k} N} \bar{\psi}_{x}) \|_{L_{t,x}^{2}} (\alpha(M) + \beta(M))
\end{equation}

\begin{equation}\label{6.108}
 \lesssim \sum_{M} (\alpha(M) + \beta(M))^{2} + C(E_{0}) \epsilon \sum_{M} (\sum_{k} \| (P_{M} \psi_{x})(P_{2^{k} M} \bar{\psi}_{x}) \|_{L_{t,x}^{2}})(\alpha(M) + \beta(M)) \lesssim \epsilon^{2}.
\end{equation}

\noindent The last inequality follows from taking $\epsilon > 0$ sufficiently small and absorbing the second term into the right hand side. This in turn implies

\begin{equation}\label{6.109}
 \sum_{N} \| (P_{M} \psi_{x})(P_{N} \bar{\psi}_{x}) \|_{L_{t,x}^{2}} \lesssim (\alpha(M) + \beta(M))^{2}
\end{equation}

\begin{equation}\label{6.110}
 + (\alpha(M) + \beta(M)) C(E_{0}) \epsilon \sum_{N} \| (P_{M} \psi_{x})(P_{N} \bar{\psi}_{x}) \|_{L_{t,x}^{2}} + C(E_{0}) \epsilon^{2} (\alpha(M) + \beta(M)).
\end{equation}

\noindent Once again absorbing the second term into the right hand side

\begin{equation}\label{6.111}
 \sum_{N} \| (P_{M} \psi_{x})(P_{N} \bar{\psi}_{x}) \|_{L_{t,x}^{2}} \lesssim \epsilon (\alpha(M) + \beta(M)).
\end{equation}

\noindent Plugging in this inequality gives the theorem. $\Box$

\section{Bilinear Estimates for $s > 0$}
Next we seek to estimate

\begin{equation}\label{7.1}
\| (P_{k} \psi_{x}(0))(P_{l} \bar{\psi}_{x}(s)) \|_{L_{t,x}^{2}}
\end{equation}

\noindent for $s > 0$. Define the double envelope at $s = 0$

\begin{equation}\label{7.2}
\gamma(t,k,l) = \sum_{j_{1}, j_{2}} 2^{-2\delta|j_{1} - k|} 2^{-2\delta|j_{2} - l|} 2^{|j_{1} - j_{2}|/2} \| (P_{j_{1}} \psi_{x}(t))(P_{j_{2}} \bar{\psi}_{x}(t)) \|_{L_{x}^{2}}.
\end{equation}

\begin{equation}\label{7.3}
 (\int \gamma(t,k,l)^{2} dt)^{1/2} \lesssim \sum_{j_{1}, j_{2}} 2^{-2\delta |j_{1} - k|} 2^{-2\delta|j_{2} - l|} 2^{|j_{1} - j_{2}|/2} \| (P_{j_{1}} \psi_{x})(P_{j_{2}} \bar{\psi}_{x}) \|_{L_{t,x}^{2}}.
\end{equation}

\noindent This implies

\begin{equation}\label{7.4}
 2^{|k - l|/2} \| (P_{k} \psi_{x})(P_{l} \bar{\psi}_{x}) \|_{L_{t,x}^{2}} \lesssim (\int \gamma(t,k,l)^{2} dt)^{1/2}.
\end{equation}

\noindent Also,

\begin{equation}\label{7.5}
 (\int \gamma(t,k,l)^{2} dt)^{1/2} \lesssim \sum_{j_{1}, j_{2}} 2^{-2\delta |j_{1} - k|} 2^{-2\delta|j_{2} - k|} \alpha(j_{1}) \alpha(j_{2}) \lesssim \alpha(k) \alpha(l).
\end{equation}

\noindent We also have the estimates

\begin{equation}\label{7.6}
 \gamma(t,k + 1,l), \gamma(t,k - 1, l) \leq 2^{2\delta} \gamma(t,k,l) \hspace{5mm} \gamma(t,k,l + 1), \gamma(t,k,l - 1) \leq 2^{2\delta} \gamma(t,k,l).
\end{equation}

\noindent Now recall Duhamel's principle. If $\psi_{x}(s)$ solves the harmonic map heat flow

\begin{equation}\label{7.7}
 \| (P_{k} \psi_{x}(s,t))(P_{l} \bar{\psi}_{x}(0,t)) \|_{L_{x}^{2}} \lesssim \| P_{k}(e^{s \Delta} \psi_{x}(0, t)) (P_{l} \bar{\psi}_{x}(0,t)) \|_{L_{x}^{2}}
\end{equation}

\begin{equation}\label{7.8}
 + \| P_{k} \partial_{l} (\int_{0}^{s} e^{(s - s') \Delta} A_{l} \psi_{x}(s', t) ds') (P_{l} \bar{\psi}_{x}(0,t)) \|_{L_{x}^{2}}
\end{equation}

\begin{equation}\label{7.9}
 + \| P_{k}  (\int_{0}^{s} e^{(s - s') \Delta} (\partial_{l} A_{l}) \psi_{x}(s', t) ds') (P_{l} \bar{\psi}_{x}(0,t)) \|_{L_{x}^{2}}
\end{equation}

\begin{equation}\label{7.10}
 + \| P_{k} (\int_{0}^{s} e^{(s - s') \Delta} A_{l} A_{l} \psi_{x}(s', t) ds') (P_{l} \bar{\psi}_{x}(0,t)) \|_{L_{x}^{2}}
\end{equation}

\begin{equation}\label{7.11}
 + \| P_{k} (\int_{0}^{s} e^{(s - s') \Delta} Im(\psi_{x} \bar{\psi_{l}})\psi_{l} (s', t) ds') (P_{l} \bar{\psi}_{x}(0,t)) \|_{L_{x}^{2}}.
\end{equation}

\begin{theorem}\label{t7.1}

 \begin{equation}\label{7.12}
\| (P_{k} \psi_{x}(s))(P_{l} \bar{\psi}_{x}(0)) \|_{L_{t, x}^{2}} \lesssim 2^{-|k - l|/2} (1 + s 2^{2k})^{-4} (\alpha(k) + \beta(k))(\alpha(l) + \beta(l)).
 \end{equation}

\end{theorem}

\noindent \emph{Proof:} We have already proved this theorem when $|k - l| \leq 10$ and for any $k,l$ when $s = 0$. As usual we start by proving

\begin{equation}\label{7.12.1}
\| (P_{k} \psi_{x}(s))(P_{l} \bar{\psi}_{x}(0)) \|_{L_{x}^{2}} \lesssim 2^{-|k - l|/2} \gamma(t,k,l).
 \end{equation}

\begin{equation}\label{7.13}
 (\ref{7.7}) \lesssim (1 + s 2^{2k})^{-4} 2^{|k - l|/2} \gamma(t,k,l).
\end{equation}

\noindent Make the bootstrap assumption

\begin{equation}\label{7.13.1}
\| (P_{k} \psi_{x}(s))(P_{l} \bar{\psi}_{x}(0)) \|_{L_{x}^{2}} \leq C \gamma(t,k,l).
\end{equation}

\begin{equation}\label{7.14}
(\int_{0}^{s} e^{-s' 2^{2k}} 2^{2k} ds')^{1/2} (\sup_{0 < s' < s} \| (P_{k - 5 \leq \cdot \leq k + 5} \psi_{x}(s'))(P_{l} \bar{\psi}_{x}(0)) \|_{L_{x}^{2}}) \| A_{x} \|_{L_{s}^{2} L_{x}^{\infty}}
\end{equation}

\begin{equation}\label{7.18}
+ (\sup_{0 < s' < s} \| (P_{k - 5 \leq \cdot \leq k + 5} \psi_{x}(s'))(P_{l} \bar{\psi}_{x}(0)) \|_{L_{x}^{2}}) \| \partial_{x} A_{x} \|_{L_{s}^{1} L_{x}^{\infty}}
\end{equation}

\begin{equation}\label{7.23}
+ (\sup_{0 < s' < s} \| (P_{k - 5 \leq \cdot \leq k + 5} \psi_{x}(s'))(P_{l} \bar{\psi}_{x}(0)) \|_{L_{x}^{2}}) \| A_{x} \|_{L_{s}^{2} L_{x}^{\infty}}^{2}
\end{equation}

\begin{equation}\label{7.27}
+ (\sup_{0 < s' < s} \| (P_{k - 5 \leq \cdot \leq k + 5} \psi_{x}(s'))(P_{l} \bar{\psi}_{x}(0)) \|_{L_{x}^{2}}) \| \psi_{x} \|_{L_{s}^{2} L_{x}^{\infty}}^{2}
\end{equation}

\begin{equation}\label{7.27.1}
\lesssim C \gamma(t,k,l) 2^{-|k - l|/2} (\| A_{x} \|_{L_{s}^{2} L_{x}^{\infty}} + \| \partial_{x} A_{x} \|_{L_{s}^{1} L_{x}^{\infty}} + \| A_{x} \|_{L_{s}^{2} L_{x}^{\infty}}^{2} + \| \psi_{x} \|_{L_{s}^{2} L_{x}^{\infty}}^{2}).
\end{equation}

\noindent Next, by Bernstein's inequality, Sobolev embedding, theorem $\ref{t6.1}$, and corollary $\ref{c6.2}$,

\begin{equation}\label{7.16}
\sum_{j \geq k + 5} 2^{k} (\sup_{0 < s' < s} \| (P_{j} \psi_{x}(s'))(P_{l} \bar{\psi}_{x}(0)) \|_{L_{x}^{2}}) \| P_{j} A_{x} \|_{L_{s}^{1} L_{x}^{\infty}}
\end{equation}

\begin{equation}\label{7.20}
+ \sum_{j \geq k + 5} 2^{k} (\sup_{0 < s' < s} \| (P_{j} \psi_{x}(s'))(P_{l} \bar{\psi}_{x}(0)) \|_{L_{x}^{2}}) \| \partial_{x} P_{j} A_{x} \|_{L_{s}^{1} L_{x}^{2}}
\end{equation}

\begin{equation}\label{7.25}
+ \sum_{j \geq k + 5}  2^{k} (\sup_{0 < s' < s} \| (P_{j} \psi_{x}(s'))(P_{l} \bar{\psi}_{x}(0)) \|_{L_{x}^{2}}) \| P_{\geq j} A_{x} \|_{L_{s,x}^{2}} \| A_{x} \|_{L_{s}^{2} L_{x}^{\infty}}
\end{equation}

\begin{equation}\label{7.29}
+ \sum_{j \geq k + 5}  2^{k} (\sup_{0 < s' < s} \| (P_{j} \psi_{x}(s'))(P_{l} \bar{\psi}_{x}(0)) \|_{L_{x}^{2}}) \| P_{\geq j} \psi_{x} \|_{L_{s,x}^{2}} \| \psi_{x} \|_{L_{s}^{2} L_{x}^{\infty}}
\end{equation}

\begin{equation}\label{7.29.1}
 \lesssim C \gamma(t,k,l) 2^{-|k - l|/2} (\| \partial_{x}^{2} A_{x} \|_{L_{s}^{1} L_{x}^{2}} + \| \partial_{x} A_{x} \|_{L_{s}^{1} L_{x}^{\infty}} + \| A_{x} \|_{L_{s}^{2} L_{x}^{\infty}} \| \partial_{x} A_{x} \|_{L_{s,x}^{2}} + \| \psi_{x} \|_{L_{s}^{2} L_{x}^{\infty}} \| \partial_{x} \psi_{x} \|_{L_{s,x}^{2}}).
\end{equation}

\noindent Also by Bernstein's inequality, Sobolev embedding, theorem $\ref{t6.1}$, and corollary $\ref{c6.2}$,

\begin{equation}\label{7.15}
\sum_{j < k - 5} \inf(2^{j} + 2^{l}, 2^{k}) (\sup_{0 < s' < s} \| (P_{j} \psi_{x}(s'))(P_{l} \bar{\psi}_{x}(0)) \|_{L_{x}^{2}}) \| P_{> k - 5} A_{x} \|_{L_{s,x}^{2}}
\end{equation}

\begin{equation}\label{7.19}
+ \sum_{j < k - 5} \inf(2^{j} + 2^{l}, 2^{k}) (\sup_{0 < s' < s} \| (P_{j} \psi_{x}(s'))(P_{l} \bar{\psi}_{x}(0)) \|_{L_{x}^{2}}) \| P_{> k - 5} \partial_{x} A_{x} \|_{L_{s}^{1} L_{x}^{2}}
\end{equation}

\begin{equation}\label{7.24}
+ \sum_{j < k - 5} \inf(2^{k}, 2^{j} + 2^{l}) (\sup_{0 < s' < s} \| (P_{j} \psi_{x}(s'))(P_{l} \bar{\psi}_{x}(0)) \|_{L_{x}^{2}}) \|  A_{x} \|_{L_{s}^{2} L_{x}^{\infty}} \| P_{\geq k - 5} A_{x} \|_{L_{s,x}^{2}}
\end{equation}

\begin{equation}\label{7.28}
+ \sum_{j < k - 5} \inf(2^{k}, 2^{j} + 2^{l}) (\sup_{0 < s' < s} \| (P_{j} \psi_{x}(s'))(P_{l} \bar{\psi}_{x}(0)) \|_{L_{x}^{2}}) \|  \psi_{x} \|_{L_{s}^{2} L_{x}^{\infty}} \| P_{\geq k - 5} \psi_{x} \|_{L_{s,x}^{2}}
\end{equation}

\begin{equation}\label{7.28.1}
 \lesssim C \gamma(t,k,l) 2^{-|k - l|/2} (\| \partial_{x} A_{x} \|_{L_{s,x}^{2}} + \| \partial_{x}^{2} A_{x} \|_{L_{s}^{1} L_{x}^{2}} + \|  A_{x} \|_{L_{s}^{2} L_{x}^{\infty}} \| \partial_{x} A_{x} \|_{L_{s,x}^{2}} + \|  \psi_{x} \|_{L_{s}^{2} L_{x}^{\infty}} \| \partial_{x} \psi_{x} \|_{L_{s,x}^{2}}).
\end{equation}

\noindent Partitioning $[0, \infty)$ and iterating over each piece proves $(\ref{7.12.1})$. To prove $(\ref{7.12})$ it remains to prove some decay in $s$ when $s > 2^{-2k}$.

\begin{equation}\label{7.31}
\aligned
\| P_{k}(\int_{0}^{(1 - \delta)s} e^{(s - s') \Delta} [2i A_{l} \partial_{l} \psi_{x}(s') - (A_{l} A_{l} + i \partial_{l} A_{l}) + Im(\psi_{x} \bar{\psi}_{l}) \psi_{l}] ds')(P_{l} \psi_{x}(0)) \|_{L_{x}^{2}} \\ \lesssim e^{-\delta s 2^{2k}} 2^{-|k - l|/2} \gamma(t,k,l) \lesssim (1 + s 2^{2k})^{-4} 2^{-|k - l|/2} \gamma(t,k,l).
\endaligned
\end{equation}

\noindent Make the bootstrap assumption

\begin{equation}\label{7.31.1}
\| P_{k}(\psi_{x}(s)) P_{l}(\psi_{x}(0)) \|_{L_{x}^{2}} \leq C(1 + s 2^{2k})^{-4} \gamma(t,k,l).
\end{equation}

\noindent When $k \leq l$,

\begin{equation}\label{7.32}
\| (P_{l} \psi_{x}(0)) (\int_{(1 - \delta)s}^{s} e^{(s - s') \Delta} P_{k}[2i \partial_{l} (A_{l} \psi_{x}) - (A_{l} A_{l} + i \partial_{l} A_{l}) \psi_{x} + Im(\bar{\psi}_{l} \psi_{x}) \psi_{l}](s') ds') \|_{L_{x}^{2}}
\end{equation}

\begin{equation}\label{7.33}
\lesssim C (1 + s 2^{2k})^{-4} \gamma(t,k,l) [\delta^{1/2} \| s^{1/2} A_{x} \|_{L_{s,x}^{\infty}} + \delta \| s \partial_{x} A_{x} \|_{L_{s,x}^{\infty}} + \delta \| s^{1/2} A_{x} \|_{L_{s,x}^{\infty}}^{2} + \delta \| s^{1/2} \psi_{x} \|_{L_{s,x}^{\infty}}^{2}]
\end{equation}

\begin{equation}\label{7.34}
\aligned
+ \sum_{j \leq k - 5} 2^{-(j - l)/2} C \gamma(t,j,l) 2^{-8k} s^{-4} [\delta^{1/2} \| s^{9/2} \partial_{x}^{8} A_{x} \|_{L_{s,x}^{\infty}} + \delta \| s^{5} \partial_{x}^{9} A_{x} \|_{L_{s,x}^{\infty}} \\ + \delta \| s^{1/2} A_{x} \|_{L_{s,x}^{\infty}} \| s^{9/2} \partial_{x}^{8} A_{x} \|_{L_{s,x}^{\infty}} + \delta \| s^{1/2} \psi_{x} \|_{L_{s,x}^{\infty}} \| s^{9/2} \partial_{x}^{8} \psi_{x} \|_{L_{s,x}^{\infty}}.
\endaligned
\end{equation}

\noindent When $k > l$,

\begin{equation}\label{7.35}
(\ref{7.32}) \lesssim C (1 + s 2^{2k})^{-4} \gamma(t,k,l) [\delta^{1/2} \| s^{1/2} A_{x} \|_{L_{s,x}^{\infty}} + \delta \| s \partial_{x} A_{x} \|_{L_{s,x}^{\infty}} + \delta \| s^{1/2} A_{x} \|_{L_{s,x}^{\infty}}^{2} + \delta \| s^{1/2} \psi_{x} \|_{L_{s,x}^{\infty}}^{2}]
\end{equation}

\begin{equation}\label{7.34}
\aligned
+ \sum_{j \leq k - 5} 2^{-|j - l|/2} C (2^{j} + 2^{l})\gamma(t,j,l) 2^{-9k} s^{-4} [\delta^{1/2} \| s^{9/2} \partial_{x}^{9} A_{x} \|_{L_{s}^{\infty} L_{x}^{2}} + \delta \| s^{5} \partial_{x}^{10} A_{x} \|_{L_{s}^{\infty} L_{x}^{2}} \\ + \delta \| s^{1/2} A_{x} \|_{L_{s,x}^{\infty}} \| s^{9/2} \partial_{x}^{9} A_{x} \|_{L_{s}^{\infty} L_{x}^{2}} + \delta \| s^{1/2} \psi_{x} \|_{L_{s,x}^{\infty}} \| s^{9/2} \partial_{x}^{9} \psi_{x} \|_{L_{s}^{\infty} L_{x}^{2}}.
\endaligned
\end{equation}

\noindent $\Box$\vspace{5mm}

\noindent We can integrate from $0$ to $s'$ with a fixed $s > 0$ in exactly the same manner. This completes the proof of theorem $\ref{t1.4}$. $\Box$

\end{document}